%\magnification =\magstep1
\baselineskip =13pt
\overfullrule =0pt
\def\op{\operatorname}

\def\K{\underline{K}}
\def\G{\underline{G}}

\def\I{\underline{I}}

\def\pf{\noindent {\sl Proof}: }

\def \s{\vskip .1cm}
\def \Z{\underline{Z}}
\def \Gr {\widehat{\op{Gr}}}
\def \T {\underline{T}}
\ifnum\pageno=1\nopagenumbers

\input amstex
\documentstyle{amsppt}
\magnification=1200
\hcorrection{.25in}
\advance\vsize-.75in
\nologo
%\noBlackBoxes
\document
 
\topmatter

\title{The elliptic curve in the S-duality theory
and Eisenstein  series for Kac-Moody groups}\endtitle
\rightheadtext{S-duality and Eisenstein series}
\author {M. Kapranov}\endauthor

\address{Department of Mathematics,
Northwestern University, Evanston, IL 60208}\endaddress

\email   kapranov\@math.nwu.edu \endemail

\abstract We establish a relation between the generating
functions appearing in the S-duality conjecture of Vafa and
Witten and geometric Eisenstein series for Kac-Moody groups.
For a pair consisting of a surface and a curve on it,
we  consider a refined generating function 
(involving $G$-bundles with parabolic structures along the curve)
which depends
on the elliptic as well as modular variables and prove
its functional equation with respect to the affine Weyl group,
thus establishing the elliptic behavior. When the curve is $\Bbb P^1$,
we calculate the Eisenstein-Kac-Moody series explicitly
and it turns out to be a certain deformation of the irreducible
Kac-Moody character, more precisely, an  analog of the Hall-
Littlewood polynomial for the affine root system. 
 We  also get an explicit formula for the universal
blowup function for any simply connected structure group.

\endabstract

%\date {???} \enddate
\endtopmatter

\document

\heading {Introduction}\endheading

\vskip .7cm

\noindent {\bf (0.1)} The goal of this paper is to develop a 
certain mathematical
framework underlying the S-duality conjecture of Vafa and Witten [VW].
 Let us recall the formulation. Let $S$ be a smooth projective surface over ${\Bbb C}$
and $G$ a semisimple algebraic group. Denote by $M_G(S, n)$
the moduli space of semistable principal $G$-bundles on $S$ with the
(second) Chern number equal to $n$ and form the generating
function of their topological Euler characteristics:
$$F_G(q) = \sum_n \chi(M_G(S,n)) q^n, \quad q=e^{2\pi i\tau}.
\leqno (1)$$
Then, the conjecture says that up to simple factors, $F_G(q)$ is a modular
form with respect to a certain 
congruence subgroup in $SL_2({\Bbb Z})$ and, moreover,
relates the image of $F_G$ under the transformation $\tau\mapsto -(1/\tau)$
with $F_{G^L}$, where $G^L$ is the Langlands dual group.
In fact, this formulation is correct only for surfaces which are not of
general type and has to be modified in general. 

Important work [Go] [LQ1-2] [Y]  has been done on the verification of this conjecture
in various particular cases (mostly restricted to $G= SL_2$ or $PGL_2$).
A treatment of the general case has been prevented by the lack
of some fundamental understanding of the problem.  More precisely,
underlying the whole theory 
(and logically preceding any finer details
of the transformation properties) is the  
 following immediate question. 

\proclaim {Question 1} Is there a  purely mathematical
 conceptual reason for $F_G$ to have
anything to do with modular functions at all? In other words,
why should the 
variable $q$ in (1), which is just a formal
variable in the generating function,  be thought of
as parametrizing  the elliptic curves $\Cal E_q= \Bbb C^*/q^{\Bbb Z}$? 
\endproclaim

Let $D = \{0<|q|<1\}$ be the punctured 
unit disk and $\Cal E\to D$ be the family of
the curves $\Cal E_q$. Our first step is to introduce a new generating
function which depends on variables living on 
(a cover of) $\Cal E$ and not just on $D$.
 Suppose $S$ is equipped with a smooth irreducible
curve $X$. We consider pairs $(P, \pi)$ where $P$ is a principal
$G$-bundle as before and $\pi$ is a parabolic structure in
$P|_X$, i.e., a reduction of the structure group from $G$
to $B$, the Borel subgroup. A $B$-bundle $Q$ on $X$ has a vector
of degrees $\op{deg}(Q)$ belonging to $L$, the lattice of coweights of $G$.
We denote by $M_{G}(n, a)$ the moduli space of semistable
parabolic $G$-bundles $P$ on $S$ with $c_2(P)=n$ and $\op{deg}(P|_X)=a$. 
Here we must  presume fixed the ``polarization parameters" [MY] neccesary to
fix the meaning of semistability. We now form the generating function
$$ E_G(q, z) =  \sum_{n\in {\Bbb Z}, a\in L} \chi(M_G(n, d)) q^n z^a,
\quad q\in D, z\in T^\vee.
\leqno (2)$$
Here $T^\vee=\op{Spec} (\Bbb C[L])$ is the ``dual torus". 
The product $D\times T^\vee$ is a natural cover of the relative
abelian variety $\Cal E_L=\Cal E\otimes_{\Bbb Z}L$, much as
 $D\times \Bbb C^*$ is
 a cover
of $\Cal E$. 

It is natural then to conjecture that in the cases when $F_G$ is
expected to have a modular behavior, $E_G$ should be, up to
simple factors, a Jacobi form [EZ], i.e., exhibit both modular behavior in
$q$ and elliptic behavior in $z$.  

\vskip .3cm

\noindent {\bf (0.2)} The first
main  result of this paper is the proof 
the elliptic behavior of the function $E_G$ in a 
wide range of ``relative situations" (allowing, for example, $G$ to
be an arbitrary simply connected group). 
 More precisely, we fix a $G$-bundle $P^\circ$
on $S-X$ and consider the moduli space $M_{G, P^\circ}(n, a)$ formed
by triples $(P, \tau, \pi)$ where $(P, \pi)$ are as above
and $\tau$ is an isomorphism of $P|_{S-X}$ with $P^\circ$.
When the self-intersection index $X^2_S=-d$ is negative, we prove 
(Theorem 2.2.1) that this space is a scheme of finite type over $k$
without imposing any stability conditions. One can form
then the relative analogs of the functions $F_G$ and $E_G$,
which we denote $F_{G, P^\circ}(q)$ and $E_{G, P^\circ}(q, z)$. 
When $X=\Bbb P^1$ and $d=1$, we have the blowup situation
considered in [LQ1-2]. 

Instead of the Euler characteristic, we work with any kind
of ``motivic measure", i.e., any invariant of algebraic
varieties additive with respect to cutting and pasting, such as, e.g.,
the number of points over the finite field. We prove (Theorem 3.4.4)
that after multiplying $E_{G, P^\circ}(q, z)$ with a certain
product over affine roots of $G$ (an analog of the Weyl-Kac denominator),
we get a series which is a regular section of the $d$th power of
a natural theta-bundle $\Theta$ on $\Cal E_L$.  This bundle is well-known
in the theory of Kac-Moody groups: among its sections one finds
 characters of level $d$
integrable representations of the Kac-Moody group associated to
the Langlands dual group
$G^L$ (note that $T^\vee$ is the maximal torus for $G^L$).

\vskip .2cm

The second main result of the paper is a complete determination of
the functions $E_{G, P^\circ}(q, z)$ and $F_{G, P^\circ}(q)$
in the case when the curve $X$ is $\Bbb P^1$. In this case,
we first show (Theorem 7.1.2) that
the type of $P^\circ$ on the punctured formal neighborhood
of $X$ can be encoded by an antidominant coweight $b$
of the Kac-Moody group of $G$, or, if one prefers, an antidominant
weight of the Kac-Moody group of $G^L$. Next we identify
(Theorem 7.3.1) the function $E_{G, P^\circ}(q, z)$
with the analog of the Hall-Littlewood polynomial [Mac]
corresponding to the affine root system of $G^L$. In other
words, this function is a certain deformation of
the character of the irreducible representation of the Kac-Moody
group with the highest weight $(-b)$. The parameter of deformation,
denoted  $\Bbb L$, is interpreted as the ``Tate motive'' (the value
of our motivic invariant on the affine line). 
In particular, we get an explicit formula for the universal
blowup function for an arbitrary simply connected structure
group. 

\vskip .3cm

\noindent {\bf (0.3)} The main idea behind our approach
 is that $E_{G, P^\circ}$
is an analog, for the Kac-Moody groups, of the unramified Eisenstein
series familiar in the Langlands' program [Lan]
 and studied in detail
for function fields by Harder [H2].

More precisely, if $X$ is a curve over
a finite field $\Bbb F_l$ and $Q$ is a principal $G$-bundle over $X$,
 then we have the associated flag fibration $\Cal F_Q$ on $X$
with fiber $G/B$. A section $s$ of $\Cal F_Q$ has a natural degree
lying in $L$, and all sections of given degree $a$ form
a scheme of finite type $\Gamma_a(Q)$.
 The Eisenstein series is just the generating
function 
$$E_{G, Q}(z) = \sum_{a\in L} \# (\Gamma_a(Q)(\Bbb F_l)) z^a, 
\quad z\in T^\vee.\leqno (3)$$
This series is known to represent a rational function which satisfies
a functional equation for each element of $W$, the Weyl group. 
When we replace $G$ by its Kac-Moody group, we get, essentially,
the series (2) and in fact Theorem 3.4.4 says that it satisfies
a functional equation for each element of $\widehat W$, the affine
Weyl group. The lattice $L$ is a subgroup of $\widehat W$,
so we get the elliptic behavior as a particular case.

Similarly, the explicit calculation of the functions $E_{G, P^\circ}$
and $F_{G, P^\circ}$ in Theorems 7.3.1 and 7.4.6 can be seen
of the analogs, in the Kac-Moody situation of Langlands'
calculation of the constant term of the Eisenstein series [Lan].

 The analogy
between $E_{G,P^\circ}$ and the ``classical" Eisenstein series
$E_{G, Q}$ is naturally understood in the context of Atiyah's work
on instantons in two and four dimensions. [A].

\vskip .3cm

\noindent {\bf (0.4)} We can now formulate succinctly our proposed
answer to Question 1: {\it  the elliptic curve appearing in the theory
of $G$-bundles on surfaces is the same curve as appears in the
theory of Kac-Moody groups}.
More precisely,  $\widehat G$,
the full Kac-Moody group
associated to $G$ (i.e., the group whose Lie algebra is associated
to the affine Cartan matrix of $G$ by the procedure of [Ka])
 is not the loop group $G((\lambda))$, nor its
central extension, but the  semidirect product
of the central extension with the group $\Bbb C^*_q$
of ``rotations
of the loop".
 Accordinly
the maximal torus is $\widehat{T}= \Bbb C^*_c\times
T\times \Bbb C^*_q$, where $\Bbb C^*_c$ is the center.
 The quotient of the open part
$\Bbb C^*_c\times
T\times D$ by $L\i \widehat W$ is a $\Bbb C^*$-bundle over $\Cal  E_L$
which corresponds to the theta-bundle, see [Lo]. So, for instance, characters
of $\widehat{G}$, being homogeneous in $\Bbb C^*_c$ and
 $\widehat W$-invariant,
are authomatically sections of a power of the theta-bundle.

In our case, however, we have a similar situation but applied to the dual
group $G^L$. Formally, the Langlands dual group to $\widehat G$
should be $\widehat{G^L}$, and 
  the center of $\widehat G$ corresponds to the loop rotation subgroup
of $\widehat{G^L}$ and vice versa. Now, the concept of the second Chern class
(for vector bundles)
is known since S. Bloch [Bl] to be intimately related to
 central extensions of matrix groups,
 and in our case this relation descends to the
central extension of $G((\lambda))$. More precisely, in our case $\lambda$
has the meaning of a local equation of $X$ in $S$, and since we consider bundles
identified with $P^\circ$ on $S-X$, the Chern class lives naturally
in the cohomology with support in $X$.

 So the variable $n$ in (1)
corresponds naturally to $\Bbb C^*_c$ in $\widehat{G}$. On the other
hand, the self-intersection index $(-d)$ of $X$ in $S$ in an invariant
of the normal bundle of $X$ in $S$ and so corresponds to $\Bbb C^*_q$
in $\widehat{G}$. When we pass to the Eisenstein series, we are transferred
to the dual group, so $d$ behaves now as if it was the central charge of
a character of $\widehat{G^L}$. In fact, out Theorem 7.3.1
identifies the Eisenstein series with the $\Bbb L$-deformation
of a character of central charge $d$.

\vskip .3cm

\noindent {\bf (0.5)} Let us now describe the contents of the paper
in more detail. In Section 1 we review the general concept
of ``motivic measures" (such as the Euler characteristic)
and of integration of such measures. Section 2 is devoted to
the finiteness theorem 2.2.1 for the moduli space $M_{G, P^\circ}(n)$.
This result can be seen as an analog of the reduction theory
developed by H. Garland [Ga] for loop groups over number fields.
Garland's twist by $e^{-tD}$, $t>0$, corresponds to our
condition that the self-intersection index of the curve is negative.
In Section 3 we formulate our results about the structure
of the generating function $E_{G, P^\circ}(q, z)$. In Section 4
we give a self-contained treatment of motivic Eisenstein series
for finite-dimensional groups. In Section 5 we
give a treatment of the Kac-Moody group as a functor
on the category of smooth  varieties. Our aim in this section
was to  make straightforward the connection with the second Chern class
which appears in the generating functions (1) and (2). 
Next, in Section 6 we 
develop the formalism of Eisenstein series for Kac-Moody groups
and prove the results of Section 3. 
Finally, Section 7 is devoted to explicit calculations of
Eisenstein-Kac-Moody series. in the case $X=\Bbb P^1$.
 We first establish a version
of Grothendieck's theorem [Gro] on $G$-bundles on $\Bbb P^1$
in the case when $G$ is replaced by the sheaf $\Cal G$ of Kac-Moody groups;
as a result, we let labelling of isomorphism classes of $\Cal G$-torsors
by antidominant affine coweights of $G$. Then, we calculate
the Eisenstein series in terms of the affine analog of the Hall
polynomial and also give a ``parahoric'' version of this calculation
which  identifies $F_{G, P^\circ}$.

\vskip .3cm

\noindent {\bf (0.6)}  This work was started 
during my stay in Max-Planck Institut f\"ur Mathematik
in 1996 and was reported at the Oberwolfach conference
on conformal field theory in June 1996. 
In writing the present
text, I was fortunate to be able to use  remarks and advice of
several people. Thus, the idea of applying the method of torus actions,
 used in Examples 2.4.2 and 3.1.2, I owe to H. Nakajima, who
uses it in a similar (but different) situation in the work in progress
with K. Yoshioka. The purely algebraic approach to the proof
of functional equations for the Eisenstein series, presented in n. 4.3,
I learned from a seminar talk of V. Drinfeld.
I am  grateful to A. Beilinson and V. Drinfeld for
numerous remarks on the previous version and to
 I. Cherednik and A. Kirillov, Jr. for discussions
of Macdonald's theory and its affine generalization. 
The revised version was prepared during my visit to MSRI in
January-February 2000. The financial support by MPI and MSRI is
hereby gratefully acknowledged.  
 The research
on this paper was partly supported by the NSF.

\vskip 2cm

\centerline {\bf \S 1. Motivic measures and zeta functions.}

\vskip 1cm

We will use a version of the formalism of ``motivic integration" from [DL].

\vskip .2cm

\noindent {\bf (1.1) Motivic measures.} 
Throughout this paper we fix a field $k$ and denote by
 $\op{Sch}_k$ the category of schemes of
 finite type over $k$. 

Let 
 $A$ a commutative ring. An $A$-valued
 {\it motivic measure} on $\op{Sch}_k$  is a function
$\mu$ which associates to any scheme $Z\in \op{Sch}_k$
 an element $\mu(Z)\in A$
depending only on the isomorphism class of $Z$ such that the two
conditions hold:

\vskip .1cm

\noindent{(1.1.1)} If $Y\i Z$ is a closed subscheme,
 then $\mu(Z)=\mu(Y)+\mu(Z-Y)$.

\vskip .1cm

\noindent{(1.1.2)} $\mu(Z_1\times Z_2) = \mu(Z_1)\times \mu(Z_2)$.

\vskip .1cm

Note that (1.1.1) implies that $\mu(Z)$ coincides with $\mu(Z_{red})$.

\vskip .2cm

\noindent {\bf (1.1.3) Examples:} (a) $k= {\Bbb C}$, $A={\Bbb Z}$ and
$\mu(V)=\chi_c(V) = \sum_i (-1)^i \dim H^i_c(V, \Bbb C)$ is the
topological Euler characteristic with compact supports. 

\vskip .1cm

\noindent (b) $k= {\Bbb F}_l$, a finite field with $l$ elements, 
$A={\Bbb Z}$, and $\mu(V) = \# V({\Bbb F}_l)$
 is the number of ${\Bbb F}_l$-points.

\vskip .1cm

\noindent (c) $k={\Bbb C}, A = {\Bbb Z}[l^{1/2}]$ is the polynomial ring,
and $\mu$ takes $V$ to its Serre 
(or virtual Hodge) polynomial $S_V(l) = \sum (-1)^i l^{i/2}
\chi (\op{gr}_i^W (H^\bullet(X(\Bbb C),
 {\Bbb C})))$. Here $W$ is the weight filtration
and $\chi$ is the Euler characteristic of the graded vector space.
If $V$ is smooth projective, then $S_V(l) =\sum_i (-1)^i l^{i/2}
\dim(H^i(V(\Bbb C), {\Bbb C}))$ is the Poincar\'e polynomial of $V$; in fact
$S_V$ is uniquely determined by this condition and the fact that it
is a motivic measure. Specializing at $l=1$, we get $S_V(1)=\chi_c(V(\Bbb C))$. 

\vskip .1cm

\noindent (d) There is the universal motivic measure $\mu_u$ whose
target $A_u$ is the ring generated by the symbols $[Z]$, where $Z$
is a quasi-projective variety, and which are subject to the relations
$[Z] = [Y] + [Z-Y]$, for a closed $Y\i Z$ and $[Z_1]\cdot[Z_2] = 
[Z_1\times Z_2]$. 

\vskip .1cm

We denote by $\Bbb L$ the value of $\mu$ at the affine line $\Bbb A^1_k$.
Thus, in the Example 1.1.3(b) $\Bbb L = l$ is the number of
the elements in $\Bbb F_l$, and in (1.1.3)(c) $\Bbb L = l$ is
the square of the generator of the polynomial ring $A=\Bbb Z[l^{1/2}]$.

\vskip .3cm

\noindent {\bf (1.2) Motivic integration.}
By $\bar k$ we denote a fixed separable closure of $k$.
Fix an $A$-valued motivic measure $\mu$. Let $Z$ be a scheme
of finite type over $k$. If $W\i Z$ is a closed subscheme
(defined over $k$), we denote by
$\bold 1_W: Z(\bar k)\to A$ the characteristic function of $Z(\bar k)$.
More generally, by a constructible ($A$-valued) function
on $Z$ we mean a function $f: Z(\bar k)\to A$ which can be represented
as a finite linear combination of characteristic
functions of closed subschemes $f=\sum_i a_i \bold 1_{W_i}$.
The set of such functions will be denoted $\op{Const}_A(Z)$.
For $f\in \op{Const}_A(Z)$ one defines its integral
 (with respect to $\mu$) as
$$\int_Z f(z) d\mu := \sum a_i \mu(W_i), \quad \op{if}\quad
f = \sum_i a_i \bold 1_{W_i}.\leqno (1.2.1)$$
It is standard  to see that this definition is independent
on the way of writing $f$ as $\sum a_i \bold 1_{W_i}$.
Further, let $\varphi: X\to Y$ be a morphism of schemes 
 and $f\in\op{Const}_A(X)$. Then the direct
image (with respect to $\mu$) of $f$ is the function
$$\varphi_*(f) = \int_{X/Y} f\cdot d\mu, \quad y\mapsto \int_{\varphi^{-1}(y)}
f(x) d\mu.\leqno (1.2.2)$$
The following is straightforward.

\proclaim{(1.2.3) Proposition} Suppose that there is a finite stratification
of $Y$ by locally closed subschemes $Y_\alpha$ such that over each $Y_\alpha$
the morphism $\varphi$ is a (Zariski) locally trivial fibration. Then
for a constructible function $f\in \op{Const}_A(X)$ the function $\varphi_*(f)$
 is also
constructible, and 
$$\int_X f(x) d\mu = \int_Y \varphi_*(f)(y) d\mu.$$

\endproclaim

\vskip .1cm

\noindent {\bf (1.3) Motivic zeta.}
 Let $X$ be a smooth algebraic variety  over $k$.
We denote by $X^{(n)}$ the $n$-fold symmetric product of $X$.
The motivic zeta-function of $X$ (associated to $\mu$) is the
formal series
$$\zeta_\mu(X, u) = \sum_{n=0}^\infty \mu(X^{(n)}) u^n \in A[[u]].\leqno 
(1.3.1)$$

\vskip .1cm

\noindent {\bf (1.3.2) Examples.} (a) If $k=\Bbb F_l$ and
$\mu$ is given by the number of $\Bbb F_l$-points, then
we get the usual Hasse-Weil zeta-function of $X$. 

(b) Consider the situation of Example 1.1.3(c) and assume that $X$
is projective. Then it is easy to see that
$$\zeta_\mu (X, u) = \prod_i (1-l^{i/2}u)^{(-1)^{i+1} b_i(X)}, \quad b_i(X)
 =\dim H^i(X, \Bbb C).
$$
This formula has the same shape as the Hasse-Weil zeta function
in the $\Bbb F_l$-case, with all the eigenvalues of the Frobenius
on $H^i$ being replaced by $l^{i/2}$.

\vskip .1cm

The following  fact generalizes the well known properties of zeta functions
of curves over $\Bbb F_l$. 

\proclaim{(1.1.9) Theorem} Let $X$ be a smooth projective
irreducible curve of genus $g$. Suppose that $A$ is a field
and $\Bbb L\neq 0$.  Suppose further that there exists a line bundle on
$X$ of degree 1. 
Then:

(a) The series $\zeta_\mu(X, u)$ represents a rational function.
In fact, $\Phi_X(u)=(1-u)(1-\Bbb Lu) \zeta_\mu(X, u)$ is a polynomial of degree $2g$. 

(b)  The function $\zeta_\mu(X,u)$ satisfies the
functional equation
$$\zeta_\mu(X, 1/\Bbb L u) = \Bbb L^{1-g} u^{2-2g} \zeta_\mu(X, u).$$

\endproclaim

\noindent {\sl Proof:} This is analogous to Artin's classical proof
for the $\Bbb F_l$-case. Let $\op{Pic}_n(X)$ be the variety of line bundles
of degree $n$ on $X$ and $p_n: X^{(n)}\to \op{Pic}_n(X)$ be the natural
projection. Clearly $p_n^{-1}(L) = \Bbb P(H^0(X, L))$, so Proposition 2.1.6
is applicable and yields
$$\zeta_\mu(X, u) = \sum_{n\in \Bbb Z} u^n \int_{L\in\op{Pic}_n(X)}
{\Bbb L^{h^0(L)-1}\over\Bbb L-1} d\mu_L.$$
This means that
$$\zeta_\mu(X, 1/\Bbb L u) = \sum_{n\in\Bbb Z} u^n \Bbb L^n \int_{L\in\op{Pic}
_{-n}(X)} {\Bbb L^{h^0(L)-1}\over\Bbb L-1}d\mu_L,$$
while
$$\Bbb L^{1-g} u^{2-2g} \zeta_\mu(X, t) = 
\sum_n u^n \int_{M\in \op{Pic}_{2g-2+n}(X)} {\Bbb L^{h^0(M)-1}\over\Bbb L -1}
\Bbb L^{1-g} d\mu.$$
Consider the isomorphism
$$\sigma: \op{Pic}_{-n}(X)\to \op{Pic}_{2g-2+n}(X), \quad L\mapsto
\omega_X\otimes L^*.$$
By the Riemann-Roch theorem, for $M=\sigma(L)$ we have
$h^0(L)-h^0(M) = -n+1-g$ and thus
$$\Bbb L^n {\Bbb L^{h^0(L)-1}\over \Bbb L-1} - \Bbb L^{1-g}
{\Bbb L^{h^0(M)-1}\over \Bbb L-1} = {\Bbb L^{1-g}-\Bbb L^n\over\Bbb L-1}.$$
So the difference between the two sides of the putative functional equation,
considered as a formal series in both positive and negative powers of $u$,
is
$${1\over \Bbb L-1} \sum_{n\in\Bbb Z} \mu(\op{Pic}_n(X))(\Bbb L^{1-g} - 
\Bbb L^n)u^n.$$
Since there exists a line bundle of degree 1 on $X$, the multiplication
by the $n$th tensor power of this bundle identifies $\op{Pic}_n(X)$
with $\op{Pic}_0(X)$ and  the above series is equal to
$${\mu(\op{Pic}_0(X))\over\Bbb L-1}\biggl(\sum_{n\in \Bbb Z} \Bbb L^{1-g} u^n -\sum_{n\in\Bbb Z}
\Bbb L^n u^n\biggr) =\mu(\op{Pic}_0(X))\biggl( {\Bbb L^{1-g}\over \Bbb L-1} \delta(u) -
{1\over\Bbb L-1}\delta(\Bbb L u)\biggr),$$
where $\delta(u) = \sum_{n\in {\Bbb Z}} u^n$ is the Fourier
series of the delta-function at 1. Now we quote the following elementary
algebraic fact, which implies our statement. 

\proclaim{(1.3.4) Lemma} Let $A$ be a field and $g_0(u)\in A((u))$,
$g_\infty(u)\in A((u^{-1}))$ be two formal Laurent series in powers
of $u$, resp. $u^{-1}$. Let $D = \sum n_i a_i$ be a positive divisor
on the affine line over $A$, with $a_i\in A, n_i\geq 0$.
 Suppose that we have an equality of formal series
in positive and negative powers of $t$:
$$g_0(u)-g_\infty(u) = \sum_i \sum_{\nu=1}^{n_i}
 c_{i\nu}\delta^{(\nu)}(u/a_i),$$
where $\delta^{(\nu)}(u)$ is the $\nu$th formal derivative of $\delta(u)$. 
Then there exists a rational function $g\in A(u)$ whose
expansion at $0$ is $g_0$, at $\infty$ is $g_\infty$
and whose divisor of poles is bounded by $D$. 
\endproclaim

\vskip .2cm

\noindent {\bf (1.3.5) Remarks.} (a) If we drop the assumption $\op{Pic}_1(X)(k)\neq\emptyset$, then the above arguments still show that
$\zeta_\mu(X,u)$ is rational and satisfies the same functional equation,
but it may have more poles. The poles will then lie at the points $u$
satisfying $u^d=1$ and $\Bbb L u^d=1$, where $d$ is such that
$\op{Pic}_d(X)(k)\neq\emptyset$. 

\vskip .1cm

(b) It is natural to expect that
the motivic zeta-functions of higher-dimensional varieties
are rational and satisfy  similar functional equations.
 
\vskip 2cm

\centerline {\bf \S 2. Relative moduli spaces of $G$-bundles.}

\vskip 1cm

\noindent {\bf (2.1) The second Chern class.} For a smooth variety $S$ over
$k$ we denote by $CH^m(S)$ the Chow group of codimension $m$ cycles on $S$
modulo rational equivalence. Thus $CH^m(S) = 
H^m(S, \K_{m, S})$, where $\K_{m, S} = K_m(\Cal O_S)$ is the sheaf
of Quillen K-functors [Q]. 

Let  $G$ be   a split simple,
simply connected
algebraic group over $k$ and $\G_S$ be the sheaf of $G$-valued regular functions
on $S$. By the work of Steinberg, Moore and Matsumoto [Ma2], as developed in
[BD] [EKLV],
  we have a natural (in $S$) central extension of sheaves 
of groups on $S$:
$$1\to \K_{2, S}\to \widetilde{\G}_S\to \G_S\to 1.\leqno (2.1.1)$$
This extension comes from a canonical element in $H^2(B_\bullet G, \K_2)$
which, in its turn, is represented by
a multiplicative $\K_{2, G}$-torsor $\Phi$ on $G$, defined uniquely
up to isomorphism of multiplicative torsors, see [BD]. In the sequel,
we will assume $\Phi$ fixed. This fixes (2.1.1) uniquely up to
a unique isomorphism.

If $P$ is a principal $G$-bundle over $S$, then (2.1.1) gives a class
$c_2(P)\in H^2(S, \K_{2, S})=CH^2(S)$ called the second Chern class of $P$.
More precisely (see [Bl]), the sheaf $\underline{P}$ of regular sections of $P$
is a sheaf of $\G_S$-torsors and the (local) liftings of
it to a sheaf of $\widetilde{\G}_S$-torsors form a $\K_{2, S}$-gerbe on
$S$, and $c_2(P)$ is the class corresponding to this gerbe by Giraud's
theory [Gi]. 

\vskip .1cm

We now describe the properties of $c_2$ which we need. 
Let $T$ be the maximal torus in $G$ and $L=\op{Hom}(\Bbb G_m, T)$
be the lattice of coweights. Let also $W$ be the Weyl group of $G$
 and $\Psi: L\times L\to \Bbb Z$ be the minimal $W$-invariant
integral negative definite scalar product on $L$. 
It is proportional to the Killing form, see [Ka], Ch. 6. 
The quadratic form $\Psi(x, x)$ is even, so $\Psi(x, x)/2$ is the
minimal $W$-invariant quadratic form on $L$. 

If $V$ is a representation of $G$ with $\dim(V)=N$, then it
induces a homomorphism $T\to\Bbb G_m^N$ of tori and at the level of
1-parameter subgroups, a homomorphism of Abelian groups
$\lambda_V: L\to\Bbb Z^N$. Let $\varphi$ be the quadratic form
on $\Bbb Z^N$ given by $\varphi(a_1, ..., a_N)=\sum_{i<j}a_ia_j$
amd $\lambda_V^*\varphi$ be the induced form on $L$. It is then an
integer multiple of the form $\Psi(x,x)/2$ and we denote by
$\varkappa_V\in\Bbb Z$ the coefficient of 
proportionality: $\varphi(\lambda_V(x)) = \varkappa_V \Psi(x, x)/2$.
We also denote by $^P V$ the vector bundle on $S$
associated to $P$ and the representation $V$. 

\proclaim {(2.1.2) Proposition} (a) The usual Chern class
$c_2({}^V P)$ is equal to $\varkappa_V c_2(P)$. 

(b) Let $Q$ be a principal $T$-bundle on $S$ and $c_1(Q)\in
H^1(S, \underline{Q})=\op{Pic}(S)\otimes_{\Bbb Z} L$ be the class
of $Q$. If $P$ is the $G$-bundle obtained from $Q$, then
$$c_2(P) = (m\otimes\Psi) (c_1(Q), c_1(Q)),$$
where $m: \op{Pic}(S)\times\op{Pic}(S)\to CH^2(S)$
is the intersection product.
\endproclaim

\noindent {\sl Proof:} Both properties follow easily from the
 construction
of the central extension in [BD], because  the form $\Psi$
is used first to construct the central extension over the torus. 
. 

\vskip .3cm

\noindent {\bf (2.2) The relative moduli space.} 
 Let now $S$ be a smooth
projective surface. We have then the degree homomorphism
$\op{deg}: CH^2(S)\to \Bbb Z$. For a principal $G$-bundle $P$ on $S$
we will abbreviate $\op{deg}(c_2(P))$ to simply $c_2(P)$ and call
it the (second) Chern number of $P$.
 Let $X\i S$ be a smooth
irreducible curve. 
We denote $X^2_S$ the self-intersection
index of $X$ in $S$ and assume that $X^2_S<0$. We denote $S^\circ = S-X$. 
Let $P^\circ$ be a principal $G$-bundle on $S^\circ$. 

\proclaim{(2.2.1) Theorem} There exists a fine moduli space
$M_{G, P^\circ}(n)$ of pairs $(P,\tau)$ where $P$ is
a principal $G$-bundle on $S$ with $c_2(P)=n$ and $\tau: P|_{S^\circ}
\to P^\circ$ is an isomorphism. This space is a scheme of finite type
 over $k$ which is empty for $n\ll 0$.
\endproclaim

\noindent {\sl Proof:} Let $\Cal M_{G, P^\circ}(n)$ be the moduli functor
on the category of $k$-schemes which corresponds to our problem.
Note that any $(P, \tau)$ as above
has no nontrivial automorphism. This means that $\Cal M_{G, P^\circ}(n)$,
coming as it does, from a stack, is in fact a sheaf with
respect to the flat topology.

\proclaim{(2.2.2) Proposition} The functor $\Cal M_{G, P^\circ}
(n)$ is represented by an ind-scheme $M_{G, P^\circ}(n)$ which is
an inductive limit of quasiprojective schemes over $k$.
\endproclaim

\noindent {\sl Proof:} First consider the case $G=SL_r$. We have then a rank
$r$ bundle $V^\circ$ on $S^\circ$ with $\det(V^\circ)\simeq \Cal O_{S^\circ}$.
Let $j: S^\circ\to S$ be the embedding. Then pairs $(P, \tau)$
as in the theorem are in bijection with locally free subsheaves $V\i
j_*V^\circ$, $\det(V^\circ)\simeq \Cal O_S$. let us assume that there
exists at least one such subsheaf, say $V_0$
(otherwise $\Cal M_{G, P^\circ}(n)$ is the functor
identically equal to the empty set and there is nothing to prove). 
Let $\bold m\i \Cal O_S$ be the sheaf of ideals of $X$.
Then any $V$ as above is contained in $\bold m^{-N}V_0$ and
contains $\bold m^NV_0$ for some $N\gg 0$. This means that
$\Cal M_{G, P^\circ}(n)$ is an inductive limit, over $N$, of
functors represented by locally closed subschemes in
$Quot(\bold m^{-N}V_0/\bold m^NV_0)$, and the case $G=SL_r$ is proved. 

The case of arbitrary $G$ is reduced to the above by taking
a sufficient number of representations $V_1, ..., V_M$ of $G$ and
realizing $M_{G, P^\circ}(n)$ as a closed sub-ind-scheme in the
product $\prod_i M_{SL_{r_i}, {}^{P^\circ}V_i}(\varkappa_{V_i} n)$. 

\vskip .15cm

Let $Y$ be the formal neighborhood of $X$ in $S$ and
$\op{Bun}_G(Y)$, $\op{Bun}_G(X)$ be the moduli stacks 
of $G$-bundles on $Y$ and $X$. We have then the restriction maps
$$\op{Bun}_{G, P^\circ}(n)\buildrel\varphi\over\longrightarrow \op{Bun}_G(Y)
\buildrel\psi\over\longrightarrow \op{Bun}_G(X).\leqno (2.2.3)$$
Here we regard $\op{Bun}_{G, P^\circ}(n)$ as a trivial stack
(a sheaf of sets). 

\proclaim{(2.2.4) Proposition} For any principal $G$-bundle $Q$
on $X$ the stack $\psi^{-1}(Q)$ is bounded
(i.e., dominated by a $k$-scheme of finite type). 
\endproclaim

\noindent {\sl Proof:} Let $X^{(d)} = \op{Spec}(\Cal O_S/\bold m^{d+1})$
be the $d$th infinitesimal neighborhood of $X$ in $S$. Given an extension
$Q^{(d)}$ of $Q$ to $X^{(d)}$, all its further extensions to
$X^{(d+1)}$, if exist, form a homogeneous space over
$H^1(X, \op{ad}(Q)\otimes (\bold m^{d+1}/\bold m^{d+2}))$.
But $\bold m^{d+1}/\bold m^{d+2} = (N^*_{X/S})^{\otimes(d+1)}$ and since
$X^2_S = \op{deg}(N_{X/S})<0$, for large $d$ the cohomology in
question vanishes. Thus all the extensions of $Q$ to a bundle
on $Y$ are determined by their restrictions to $X^{(d)}$ for some $d$,
and the latter form a bounded family.

\vskip .1cm

\proclaim{(2.2.5) Proposition} For any principal $F$-bundle $\widehat P$
on $Y$, the stack $\varphi^{-1}(\widehat P)$ is bounded. 
\endproclaim

\noindent {\sl Proof:} Let $Y^\circ = Y\cap S^\circ$ be the punctured
formal neighborhood of $X$ in $S$. Supposing $\varphi^{-1}(\widehat P)\neq
\emptyset$ (otherwise there is nothing to prove), there is a $G$-bundle
$P$ on $S$ whose restriction to $Y$ is isomorphic to $\widehat P$ and whose
restriction to $S^\circ$ is isomorphic to $P^\circ$. Any other object
of $\varphi^{-1}(\widehat P)$ is then obtained by gluing $\widehat P$ and $\widehat P^\circ$
over $Y^\circ$ in a different way. This new gluing function is an
element of $\op{Aut}(P|_{Y^\circ})$.
If $g, g'\in \op{Aut}(P|_{Y^\circ})$ are such that $g'=gh$ with
$h\in\op{Aut}(\widehat P)$, then the gluings by $g$ and $g'$ give the same
object (element) of $M_{G, P^\circ}(n)$. Thus it is enough to
prove the following.

\proclaim{(2.2.6) Proposition} Let $P$ be any principal $G$-bundle
on $S$. Then the quotient $\op{Aut}(P|_{Y^\circ})/\op{Aut}(P|_Y)$
is bounded.
\endproclaim

\noindent {\sl Proof}: We first consider the case $G=SL_r$, so $P$ is
given by a rank $r$ bundle $V$ on $S$. As $X^2_S<0$, the curve $X$
can be blown down to a (possibly singular) point $p$ on a new surface
$S'$. Let $\sigma: S\to S'$ be the projection. Because $\sigma^{-1}$
identifies $S'-\{p\}$ with $S^\circ = S-X$, we have the
embeddings $j: S^\circ\to S$ and $j'=\sigma j: S^\circ\to S'$. 
Let $Y'$ be the formal neighborhood of $p$ in $S'$ and
$\widehat\sigma: Y\to Y'$ be the projection. Let also
${Y'}^{\circ} = Y'-\{p\}$ be the punctured formal neighborhood.
Denote $\bold m\i \Cal O_S$ and $\widehat {\bold m}\i\Cal O_Y$ the sheaves
of ideals corresponding to $X$. The morphism $\sigma$ is projective,
in fact
$$S=\op{Proj}\biggl(\bigoplus_{d\geq 0} \sigma_*(\bold m^d)\biggr)$$
and the relative sheaf $\Cal O(1)$ on the projective spectrum
is identified with $\bold m$. The same is true for $\widehat\sigma$.
Now, let us aply the relative version of Serre's theorem 
(about the equivalence of the categories of coherent sheaves
and graded modules)
to the coherent sheaf $\widehat V=V|_Y$ on $Y$. We conclude that
for $d\gg 0$ we have
$$\op{Hom}_Y(\widehat V, \widehat V) = \op{Hom}_{Y'}\bigl(
\widehat \sigma_*(\widehat V(d)), \widehat \sigma_*(\widehat V(d))\bigr),$$
where $\widehat V(d) = \widehat V\otimes\Cal O(d) = \bold m^d \widehat V$.
Hence $\op{Aut}(\widehat V) = \op{Aut}(\widehat \sigma_*\widehat V(d))$. 
Since $\op{Aut}(\widehat V) =
\op{Aut}(\widehat V(d))$, we can replace $\widehat V$ by $\widehat V(d)$ and assume that
$\op{Aut}(\widehat V) = \op{Aut}(\widehat\sigma_*(\widehat V))$. On the other hand,
$$\op{Aut}(\widehat V|_{Y^\circ}) = \op{Aut}( (\widehat\sigma_*\widehat V)|_{Y'\circ})
\i \op{Aut}((j'_*{j'}^*\sigma_*V)|_{Y'}).$$
Since $j'$ is an embedding of the complement of a point into a surface,
$j'_*{j'}^*\sigma_*V$ is a coherent sheaf on $S'$ which contains the
(torsion-free) sheaf $\sigma_*V$ and coincides with it outside $p$.
We have therefore an exact sequence
$$0\to \widehat\sigma_*\widehat V\to (j'_*{j'}^*\sigma_*V)|_{Y'}\to\Cal F\to 0,
\leqno (2.2.7)$$
where $\Cal F$ is a coherent sheaf supported at $p$.  
It follows that the quotient of the automorphism group of the
middle sheaf by the automorphism group of the left sheaf is bounded. 
Thus Proposition 2.2.6 is proved for $G=SL_r$. The case of general $G$ is deduced from this by applying the above reasoning to vector bundles
associated to sufficiently many representations of $G$. 

\vskip .3cm

\noindent {\bf (2.3) End of the proof of Theorem 2.2.1.}
because of Propositions 2.2.4 and 2.2.5, we are reduced to the
following fact. 

\proclaim{(2.3.1) Proposition} The image of $M_{G, P^\circ}(n)$
under the map $\psi\varphi$ in (2.2.3) is bounded.
\endproclaim

\pf As before, we start with the case $G=SL_r$, so we work with rank $r$ bundles
$V$ on $S$ with $\det(V)\simeq \Cal O$, $c_2(V)=n$ equipped with an identification $V|_{S^\circ}\to V^\circ$ where $V^\circ$ is some fixed
bundle on $S^\circ$. We will use the notation introduced in the proof
of Proposition 2.2.6, in particular, we will use the global analog
of the sequence (2.2.7), which we write in the form
$$0\to\sigma_*V\to j'_*V^\circ \to\Cal F\to 0\leqno (2.3.2)$$
with $\Cal F$ supported at $p\in S'$. We choose an embedding
$f: S'\to \Bbb P^N$ and apply the Grothendieck-Riemann-Roch theorem
to $V$ and to $f\sigma: S\to \Bbb P^N$. From (2.3.2) we infer that
$$ch_N(f_*\sigma_*V) = c-l(\Cal F),$$
where $l(\Cal F)=\dim\, H^0(\Cal F)$ is the length of the 0-dimensional
sheaf $\Cal F$ and $c$ is a constant depending only on $V^\circ$
but not
on $V$. The only higher direct image to take into account is
$R^1(f\sigma)_*V= f_* R^1\sigma_*V$. The sheaf $R^1\sigma_*V$ is again
supported at $p$ and we have
$$ch_N(f_*R^1\sigma_*V) = l(R^1\sigma_*V).$$
As $c_1(V)=0, c_2(V)=n$ are fixed, the GRR theorem
implies part (a) of the following fact.

\proclaim{(2.3.3) Lemma}(a)  For $V\in M_{SL_r, V^\circ}(n)$ the sum
$l(\Cal F)+l(R^1\sigma_*V)$ is equal to some fixed constant
$c=c(V^\circ, n)$. In particular,
$l(R^1\sigma_*V)$ is bounded. 

(b) Further, for $n\ll 0$ we have
$c(V^\circ, n)<0$ which implies that $M_{SL_r, V^\circ}(n) = 
\emptyset$.
\endproclaim

Part (b) is true because the dependence of
$c(V^\circ, n)$ on $n$ comes out to be  affine linear
with a positive coefficient in the linear part. 

\vskip .1cm

Now, since $R^1\sigma_*V$ is supported at $p=\sigma(X)$, we have
$$H^0(R^1\sigma_*V) = H^1(Y, V|_Y),\leqno (2.3.4)$$
where $Y$ is, as before, the formal neighborhood of $X$ in $S$. In fact, one
can replace in (2.3.4) $Y$ by the infinitesimal neighborhood $X^{(d)}$
for sufficiently large $d$. Consider the spectral sequence corresponding
to the filtration of $V|_Y$ (or $V|_{X^{(d)}}$) by powers of $\bold m$.
Its quotients are $V|_X\otimes (N^*_{X/S})^{\otimes i}$, $i\geq 0$.
By analyzing the spectral sequence, we find that the component $H^1(X, V|_X)$
of the $E_1$-term consists of permanent cycles, by dimension reasons. So
$\dim \, H^1(X, V|_X)\leq l(R^1\sigma_*V)\leq c(V^\circ, n)$ bounded.
As $V|_X$ has $c_1=0$, we find, by Kleiman's criterion [Kl]
that all possible $V|_X$ for $V\in M_{SL_r, V^\circ}(n)$
form a bounded family. 

This finishes the proof for $G=SL_r$. The case of arbitrary $G$ is reduced
to this one by considering associated vector bundles and applying
the following fact which is an easy consequence of the reduction theory
of Harder [H1]

\proclaim{(2.3.5) Lemma} A family of principal $G$-bundles on a curve is bounded
if and only if for any representation of $G$ the corresponding family of
vector bundles is bounded.
\endproclaim

\vskip .2cm

\noindent {\bf (2.4) The generating function.}
We now fix a ring $A$ and an $A$-valued motivic measure $\mu$ on $\op{Sch}_k$,
as in \S 1. To the data of $S, X, G$ and $P^\circ$ we can, in virtue of
Theorem 2.2.1, associate the generating function
$$F(q) = F_{G, P^\circ}(q) = \sum_{n\in\Bbb Z} \mu(M_{G, P^\circ}(n))q^n
\quad\in\quad A((q)).\leqno (2.4.1)$$

\noindent {\bf (2.4.2) Example.} Let $k=\Bbb C$, $A=\Bbb Z$ and $\mu$
be given by the Euler characteristic, as in Example 1.1.3(a). Assume that
$P^\circ$ is trivial. Then we have a $G$-action on $M_{G, P^\circ}(n)$
given by 
$g(P, \tau) = (P, g\circ\tau)$. Let us consider the induced $T$-action and
find the fixed locus $M_{G, P^\circ}(n)^T$. By definition, it consists of
$(P, \tau)$ which are isomorphic (as bundles with trivialization)
to $(P, t\circ\tau)$ for any $t\in T$.
Since the moduli problem associated to $M_{G, P^\circ}(n)$
has trivial automorphism groups, we find that $M_{G, P^\circ}(n)^T$ consists of
$(P,\tau)$ which come from $T$-bundles on $S$ trivialized on $S^\circ$.
In the identification $\op{Bun}_T(S)=\op{Pic}(S)\otimes_{\Bbb Z}L$,
bundles trivial on $S^\circ$ form a subgroup $\bold m\otimes L$
where $\bold m$ is, as before, the (invertible) sheaf of ideals of $X$.
We will denote such bundles by $\bold m\otimes a$, $a\in L$. The trivialization
of $\bold m\otimes a$ over $S^\circ$ is unique, up to isomorphism in our sense,
and the second Chern number of the associated $G$-bundle is, by 
Proposition 2.1.2(b), equal to $(-d)\Psi(a,a)$, where $d=-X^2_S$. Therefore
$$\chi(M_{G, P^\circ}(n)) = \chi(M_{G, P^\circ}(n)^T) = 
\#\{a\in L: -d\Psi(a,a)=m\},$$
and hence
$$F_{G, P^\circ}(q) = \sum_{a\in L}q^{-d\Psi(a,a)/2} = \theta_L(q^d)$$
is the theta-function (or, rather the theta-zero-value) of the
lattice $L$ with the positive definite quadratic form $(-\Psi(x,x)/2)$. 

\vskip 2cm

\centerline {\bf \S 3. The refined generating function and its elliptic behavior.}

\vskip 1cm

\noindent {\bf (3.1) The relative moduli space of parabolic bundles.}
We keep the notation of \S 2. A $T$-bundle $U$ on $X$ has a 
degree $\op{deg}(U)\in L$ which is the image of the class of $U$ in
$\op{Bun}_T(X) = \op{Pic}(X)\otimes L$ under $\op{deg}' \otimes \op{Id}$,
where $\op{deg}': \op{Pic}(X)\to \Bbb Z$ is the usual degree of line bundles.
Let $B$ be a fixed Borel subgroup in $G$ containing $T$. As $B/[B,B]\simeq T$,
a $B$-bundle $R$ gives a $T$-bundle, whose degree will also  be denoted
$\op{deg}(R)\in L$. We set $d=-X^2_S$. By our assumptions, $d>0$. 
Let $L_+\i L$ be the semigroup spanned by the positive coroots. 
For a principal $G$-bundle $Q$ on $X$ and $a\in L$ let $\Gamma_a(Q)$ be the scheme 
of all $B$-structures in $Q$ of degree $a$. The following result is due to G.
Harder [H2]. 

\proclaim {(3.1.1) Proposition}
 $\Gamma_a(Q)$ is a quasiprojective scheme over $k$. There is $a_0
in L$ such that $\Gamma_q(Q)=\emptyset$ unless $a\in a_0+L_+$. 

\endproclaim

Together with Theorem 2.2.1 this implies the following. 

\proclaim {(3.1.2) Corollary}
 For any $a\in L, n\in\Bbb Z$ there exists a $k$-scheme $M_{G, P^\circ}(n, a)$
of finite type, which is a fine moduli space of triples $(P, \tau, \pi)$ where:

$\bullet$ $P$ is a $G$-bundle on $S$ with $c_2(P)=n$;

$\bullet$ $\tau: P|_{S^\circ}\to P^\circ$ is an isomorphism;

$\bullet$ $\pi$ is a $B$-structure in $P|_X$ of degree $a$.

\endproclaim

 We now  form the generating function
$$E_{G, P^\circ}(q, z, v) = \sum_{n\in \Bbb Z, \, a\in L} \mu(M_{G, P^\circ}
(n, a))q^n z^a v^d.\leqno (3.1.3)$$
Here $\mu$ is a fixed $A$-valued motivic measure on $\op{Sch}_k$
and $z$ is a formal variable running in the ``dual torus"
$T^\vee = \op{Spec} \, \Bbb Z[L]$ whose lattice of characters is $L$. 
The variable $v$ (in which $E_{G, P^\circ}$ is homogeneous of
degree $d$) is added for convenience of future formulations. 

\vskip .1cm

\noindent {\bf (3.1.4) Example.} Consider the situation of Example 2.4.2.
As before, we have a $G$-action on $M_{G, P^\circ}(n, a)$. The set 
of $T$-fixed points consists  of one element, the $G$-bundle corresponding
to $\bold m\otimes a$, if $\Psi(a,a) = n$, and is empty otherwise. This
implies that
$$E_{G, P^\circ}(q, z, v) = \sum_{a\in L} q^{-d\Psi(a,a)/2} z^a v^d = \theta_L(q^d, z) v^d,$$
where $\theta_L(q, z)$ is now the full theta-function, depending on the 
elliptic variables as well as on the modular one. 

\vskip .3cm

\noindent {\bf (3.2) The affine root system.} We denote by $L^\vee=
\op{Hom}(L, \Bbb Z)$ the weight lattice of $G$ and
by $\Delta_{\op{sim}}\i \Delta_+ \i\Delta\i L^\vee$  the sets of simple,
resp. positive, resp. all roots of $G$. Let
 $\rho\in L^\vee$ is the half-sum of the positive roots.
Let   $L^\vee_{\op{aff}} = {\Bbb Z}\oplus L^\vee\oplus \Bbb Z$ be
the lattice of affine weights  and 
$$\widehat \Delta = \{ (0, \alpha, n), \alpha\in \Delta, n\in {\Bbb Z}\} \i
 L^\vee_{\op{ aff}},$$
$$\widehat  \Delta_+ = (\{0\}\times \Delta_+\times \{0\}) \,\, \cup \,\, (\{0\}\times
 \Delta\times \Bbb Z_{\geq 0})$$
be the system of affine roots of $G$, resp. positive affine roots. The
set of simple affine roots will be denoted by
$$\widehat  \Delta_{\op{ sim}} = (\{0\}\times \Delta_{\op{ sim}}\times\{0\})
 \,\,\cup\,\, \{(0, -\theta, 1,$$
where $\theta\in\Delta_+$ is the maximal root.
We denote $\widehat\rho = (0,\rho, h^\vee)\in L^\vee_{\op{aff}}$,
where $h^\vee$ is the dual Coxeter number of $G$, see [Ka].

\vskip .1cm

As $G$ is assumed simply connected, its coweight lattice $L$ coincides
with the coroot lattice. For any $\alpha\in\Delta$ let $\alpha^\vee\in L$
be the corresponding coroot and $\rho^\vee=(1/2)\sum_{\alpha\in\Delta_+}
\alpha\in (1/2)L$. We set $L_{\op{aff}}=\Bbb Z\oplus L\oplus \Bbb Z$.
This is the lattice dual to $L^\vee_{\op{aff}}$.
Denote by $\widehat\Delta^\vee\i L_{\op{aff}}$ the system of affine
coroots. Again, we use the notation $\alpha^\vee$ to denote the affine
coroot corresponding to $\alpha\in\widehat\Delta$.
Thus, if $\alpha =(0, \beta, n)$ with $\beta\in\Delta$ and $n\in\Bbb Z$,
then $\alpha^\vee = (n, \beta^\vee, 0)$. 
 We also write
$\widehat\rho^\vee = (h, \rho^\vee, 0)$, where
$h$ is the Coxeter number of $G$. 

\vskip .3cm

\noindent {\bf (3.3) The abelian variety $\Cal E_L$.}
We now recall the standard appearance of elliptic curves in
the Kac-Moody theory, see [Lo].

 Let 
$\widehat  W_{} = W\ltimes L$ be the affine Weyl group of $G$.
It is generated by the
reflections $s_\alpha, \alpha\in\widehat\Delta_{\op{sim}}$. 
This group acts on $L_{\op{aff}} = \Bbb Z\oplus L\oplus \Bbb Z$ by 
$$w\circ (n,a, m) = (n, w(a), m),\quad  w\in W, $$
$$b\circ (n,a, m) =
 \biggl(n+ \Psi(a,b)+ {1\over 2} \Psi(b, b)m ,  a+mb, m\biggr),
\quad  b\in L, \leqno (3.3.1)$$
see [PS]. Note that this action preserves the
 subgroup $\Bbb Z\oplus L\i L_{\op{aff}}$
given by $m=0$. Accordingly, we have the
 $\widehat W$-action on $T^\vee_{\op{aff}} = 
\Bbb G_m\times T^\vee\times \Bbb G_m = \op{Spec}\, \Bbb Z[L_{\op{aff}}]$.
We will denote a typical point of $T^\vee_{\op{aff}}$
 by $t=(q, z, v)$, where $q, v\in\Bbb G_m$
and $z\in T^\vee$. because the sugroup $\{m=0\}$ is
 preserved by $\widehat W$, we have a 
$\widehat W$-action on the torus $\Bbb G_m\times T^\vee = \op{Spec}\, \Bbb Z[\Bbb Z\oplus
L]$ with coordinates $q, z$, so that the projection $T^\vee_{\op{aff}}\to
\Bbb G_m\times T^\vee$ given by $(q,z,v)\mapsto (q,z)$, is $\widehat W$-equivariant.
In other words, $T^\vee_{\op{aff}}$ is a $\widehat W$-equivariant $\Bbb G_m$-bundle
on $\Bbb G_m\times T^\vee$. We denote by $\theta$ the corresponding $\widehat W$-equivariant
line bundle on $\Bbb G_m\times T^\vee$ (so sections of $\theta$  are functions $f(q,z,v)$
homogeneous in $v$ of degree 1).

Assume now that the ring $A$ where our motivic measure $\mu$ takes values,
is a field. Then $A((q))$ is a complete discrete valued field.
Let $\Cal E$ be the Tate elliptic curve over $A((q))$.
Thus $\Cal E^{\op{an}}$, the rigid analytic $A((q))$-space corresponding
to $\Cal E$, is the quotient $\Cal E^{\op{an}} = \Bbb G_{m, A((q))}^{\op{an}}
/q^{\Bbb Z}$. Consider also the abelian variety $\Cal E_L$
over $A((q))$ such that
$$\Cal E_L^{\op{an}} = T^{\vee, \op{an}}_{A((q))}/q^{\Psi(L)}.\leqno
(3.3.2)$$
Here we regard $\Psi$ as a homomorphism $L\to L^\vee$ and view $L^\vee$
as the lattice of 1-parameter subgroups in $T^\vee$, so
$q^\lambda$, $\lambda\in L^\vee$, is the value at $q$ of the
1-parameter subgroup $\lambda$. Note that $T^\vee_{A((q))}$ can be viewed
as a completion of $(\Bbb G_m\times T^\vee)_A$, with $q$ being the coordinate in $\Bbb G_m$.
The action of $L\i \widehat W$ on $T^\vee_{A((q))}$
 coming, via  this identification, from the action on $\Bbb G_m\times T^\vee$
induced by (3.3.1), is the same as the one used in (3.3.2).
In particular, $W$ acts on $\Cal E_L$ and
$$(\Cal E_L/W)^{\op{an}} = T^{\vee, \op{an}}_{A((q))}/\widehat W. \leqno
(3.3.3)$$
The $\widehat W$-equivariant line bundle $\theta$ on $\Bbb G_m\times T^\vee$
gives a $W$-equivariant line bundle on $\Cal E_L$, which we denote by $\Theta$.
Similarly, 1-cocycles of $\widehat W$ with coefficients in $A[T^\vee_{\op{aff}}]^*$,
the group of invertible regular functions, give rise to 
line bundles on $\Cal E_L/W$ and
hence to $W$-equivariant line bundles on $\Cal E_L$. We will need the following
two cocycles:

\vskip .2cm

\noindent {(3.3.4)} The cocycle $w\mapsto (-1)^{l(w)}$, where $l(w)$ is the
length of $w$. This cocycle is trivial on $L\i \widehat L$.

\vskip .1cm

\noindent (3.3.5) The cocycle
$$w\mapsto t^{w(\widehat\rho^\vee)-\widehat\rho^\vee}=
\prod_{\alpha\in\widehat\Delta_+ \cap  w^{-1}(\widehat \Delta_-)}
t^{\alpha^\vee}, \quad t = (q, z, v)\in T^\vee_{\op{aff}}.$$
As $2\widehat\rho^\vee\in L_{\op{aff}}$, the square of this cocycle is
a coboundary. We denote by $\Cal L$ the line bundle on $\Cal E_L$
corresponding to this cocycle. Thus $\Cal L^{\otimes 2}\simeq \Cal O$.

\vskip .3cm

\noindent {\bf (3.4) The functional equation.} We now proceed to formulate
the first main result of this paper. We write $\zeta(u)$ for $\zeta_\mu(X,u)$,
the motivic zeta function of $X$ and assume from now on that
$\op{Pic}_1(X)(k)\neq\emptyset$. This condition, added for simplicity
of formulations, always holds when $k$ is algebraically
closed or finite. Thus $\Phi(u) = \Phi_X(u)$, the numerator of
$\zeta(u)$, is a polynomial if degree $2g$ with constant term 1.
As before, we denote $d=-X^2_S$ and assume $d>0$. We also
assume that $A$ is a field and $\Bbb L\neq 0$.

\proclaim {(3.4.1) Theorem} The series $E_{G, P^\circ}(t) = 
E_{G, P^\circ}(q, z, v)$ converges to a meromorphic function on
$(T^\vee\times \Bbb G_m)^{ \op{an}}_{A((q))}$
(homogeneous of degree $d$ in the variable $v$) which satisfies, for any $w\in\widehat W$, the 
functional equation
$$E_{G, P^\circ}(t) = \Lambda_w(t) E_{G, P^\circ}(\Bbb L^{\widehat\rho -
w(\widehat
\rho)} t),$$
where
$$\Lambda_w(t) = \Bbb L^{l(w)(1-g)}\prod_{\alpha\in\widehat\Delta_+\cap
w^{-1}(\widehat\Delta_-)}{\zeta(\Bbb Lt^{\alpha^\vee})\over
\zeta(t^{\alpha^\vee})}.$$
Further, the singularities of $E_{G, P^\circ}(t)$ are contained
among the singularities of the rational functions $\Lambda_w(t)$. 
\endproclaim

It is convenient to reformulate this theorem as follows. Let
$$D(t) = \prod_{\alpha\in\widehat\Delta_+} (1-\Bbb L ^2t^{\alpha^\vee})
\Phi(t^{\alpha^\vee}).\leqno (3.4.2)$$
As with the Weyl-Kac denominator [Ka], it is clear that
$D(t)$ is an analytic function on $T^{\vee, \op{an}}_{A((q))}$
(it does not depend on $v$),  a kind of
theta-function. Then set
$$N_{G, P^\circ}(t) = E_{G, P^\circ}(t) D(t)\leqno (3.4.3)$$
and call $N_{G, P^\circ}(t)$ the
{\it numerator} of $E_{G, P^\circ}$. Taking into account the functional
equation for $\zeta(u)$ together with the fact that $D(t)$ dominates all
the poles of all the $\Lambda_w(t)$, we can reformulate Theorem 3.4.1 as
follows.

\proclaim{(3.4.4) Theorem} The series $N_{G, P^\circ}(t)$ is well-defined
and represents an analytic function on $T^{\vee, \op{an}}_{A((q))}$ which,
for any $w\in\widehat W$, satisfies the functional equation
$$N_{G, P^\circ}(t) = (-1)^{l(w)}
\bigl( (\Bbb L^{\widehat\rho} t)^{w(\widehat\rho^\vee)-
\widehat\rho^\vee}\bigr)^{1+2g} N_{G, P^\circ}(\Bbb L^{\widehat\rho-w(\widehat\rho)} w(t)).$$
\endproclaim

In other words, up to a shift and a multiplication by a monomial,
$N_{G, P^\circ}$ is a regular section of the line bundle $\Cal L\otimes\Theta^{
\otimes d}$ on $\Cal E_L$, which is antisymmetric with respect
 to the $\widehat W$-action on this equivariant bundle.

This identifies $N_{G, P^\circ}$ (and therefore $E_{G, P^\circ}$) as
an element of a finite-dimensional $A((q))$-vector space depending only
on $G$ and $d$. 

The proofs of Theorems 3.4.1 and 3.4.4 will be given in \S 6. 

\vskip 2cm

\centerline {\bf \S 4. Motivic Eisenstein series for $G$.}

\vskip 1cm

We first summarize the classical theory of unramified geometric
Eisenstein series [H2]
putting it into the motivic framework. In this section $G$ will be assumed to
be any split reductive group over $k$ such that $[G, G]$
is simply connected. The notations pertaining to the root system remain the
same.

\vskip .3cm

\noindent {\bf (4.1) Rationality and functional equation.} Let 
$Q$ be a principal
$G$-bundle on $X$.
The motivic Eisenstein series associated to $Q$, is the generating
function
$$E_{G, Q}(z) = \sum_{a\in L} \mu(\Gamma_a(Q)) z^a, \quad z\in T^\vee,
\leqno (4.1.1)$$
where $\Gamma_a(Q)$ is the scheme of $B$-structures on $Q$ of degree $a$,
see (3.1) for this and other notations.

\proclaim{(4.1.2) Theorem} (a) The series $E_{G, Q}(z)$ represents
a rational function on $T^\vee_A$.

(b) This rational function satisfies, for any $w\in W$, the functional
equation
$$E_{G, Q}(z) =   M_w(z) E_{G, Q}(\Bbb L^\rho w(\Bbb L^{-\rho}z)),$$
where
$$M_w(z) =  \Bbb L^{l(w)(1-g)}\prod_{\alpha\in\Delta_+\atop
w(\alpha)\in\Delta_-}{\zeta(\Bbb Lz^{\alpha^\vee})\over\zeta(z^{\alpha
^\vee})}.$$
Further, the singularities of $E_{G, Q}(z)$ are contained
among the singularities of the rational functions $M_w(t)$, $w\in W$.
\endproclaim

Let us give an equivalent formulation in terms of the ``numerator"
$$N_{G, Q}(z) = E_{G, Q}(z) \cdot\prod_{\alpha\in\Delta_+} 
(1-\Bbb L ^2z^{\alpha^\vee})
\Phi(z^{\alpha^\vee}).\leqno (4.1.3)$$

\proclaim{(4.1.4) Theorem} The series $N_{G, Q}(z)$ is a Laurent
polynomial which satisfies, for any $w\in W$, the functional equation
$$N_{G, Q}(z) = (-1)^{l(w)} 
 \bigl( (\Bbb L^\rho z)^{w(\rho^\vee)-
\rho^\vee}\bigr)^{1+2g} N_{G, Q}(\Bbb L^{\rho-w(\rho)} w(z)).$$
\endproclaim

The equivalence of (4.1.2) and (4.1.4) is a formal consequence of the
functional equation for $\zeta(u)$ and is left to the reader.

\vskip .3cm

\noindent {\bf (4.2) Proof for $G=GL_2$.} Assume $G=GL_2$. Thus $Q$ comes from 
a rank 2 vector bundle $V$ on $X$. In this case $L=\Bbb Z^2$ and
a $B$-structure on $Q$ of degree $a=(a_1, a_2)$ is the same as a rank
1 subbundle $\Cal M\i V$ with $\op{deg}(\Cal M)=a_1$ and
$\op{deg}(V/\Cal M)=a_2$. Thus $\Gamma_a(Q)$ parametrizes such subbundles.
Recall that a subbundle is just a coherent subsheaf which is locally
a direct summand. Any proper coherent subsheaf $\Cal M\i V$
(subbundle
or not) is automatically locally free of rank 1, so $\op{deg}(\Cal M)$
is defined. Further, $V/\Cal M$ is a direct sum $(V/\Cal M)_{\op{lf}}\oplus
(V/\Cal M)_{\op{tors}}$ of a locally free sheaf of rank 1 and a
torsion sheaf. We define
$$\op{deg}(V/\Cal M) = \op{deg}((V/\Cal M)_{\op{lf}})+ \dim\, H^0((V/\Cal M)_{
\op{tors}}).$$
Let $\overline{\Gamma_{a_1, a_2}}(Q)$ be the scheme parametrizing
subsheaves $\Cal M\i V$ with $\op{deg}(M)=a_1$ and
$\op{deg}(V/\Cal M) = a_2$. This is just a component of the
$Quot$ scheme of $V$.

As $L=\Bbb Z^2$, we have $T^\vee = \Bbb G_m^2$ and write
$E_{G, Q}(z)$ as $E_{G, Q}(z_1, z_2)$. 

\proclaim {(4.2.1) Proposition} (cf. $\op{[Lau]}$.) The generating function
$$\widetilde E(z_1, z_2) = \sum_{a_1, a_2\in\Bbb Z} \mu (\overline{\Gamma_{a_1, a_2}}(Q)) z_1^{a_1} z_2^{a_2}$$
is equal to $\zeta(z_2/z_1) E_{G, Q}(z_1, z_2)$.
\endproclaim 
\pf Any rank 1 subsheaf $\Cal M\i V$ can be uniquely
written in the form the form $\Cal M=\overline{\Cal M}
(-D)$ where $\overline{\Cal M}$ is a rank 1 subbundle and 
$D$ is a positive divisor on $X$. This implies the statement.

\vskip .1cm

In view of the functional equation for $\zeta(u)$, Theorem 4.1.2 for
the $GL_2$-bundle $Q$ is equivalent to the following statement. 

\proclaim{(4.2.2) Proposition} The series $\widetilde E(z_1, z_2)$
represents a rational function
with only poles being simple poles along the lines $z_1=z_2$ and
$z_1=\Bbb L^2 z_2$. It satisfies the functional equation
$$\widetilde{E}(z_1, z_2) = (\Bbb L z_1/z_2)^{2-2g} \widetilde E(\Bbb L z_2, \Bbb L^{-1}z_1).$$
\endproclaim

\pf We will use Lemma 1.3.4 and prove that the difference of the
two sides of the proposed equality, considered as a formal
series in $z_i, z_i^{-1}$, is a sum of two delta functions. 
Consider the projection
$$p_{a_1, a_2}: \overline{\Gamma}_{a_1, a_2}(Q)\to\op{Pic}_{a_1}(X),$$
which  takes $\Cal M\i V$ into the isomorphism class of $\Cal M$ in
the Picard group.  If $a_1+a_2=\op{deg}(V)$
(otherwise $\overline{\Gamma}_{a_1, a_2}(Q)=\emptyset$),
then the fiber $p_{a_1, a_2}^{-1}(\Cal M)$
is the projective space $\Bbb P(\op{Hom}(\Cal M, V))$. 
Thus the coefficient at $z_1^{a_1}z_2^{a_2}$, $a_1+a_2=\op{deg}(V)$ in
$\widetilde{E}(z_1, z_2)$ is
$$\int_{\Cal M\in \op{Pic}_{a_1}(X)} {\Bbb L^{\dim\,\op{Hom}(\Cal M, V)}-1
\over \Bbb L -1} d\mu_{\Cal M},$$
while the coefficient at the same monomial in
$(\Bbb Lz_1/z_2)^{2-2g} \widetilde E(\Bbb Lz_2, \Bbb L^{-1} z_1)$ is
 $$\int_{\Cal M'\in \op{Pic}_{a_2+2-2g}(X)} \Bbb L^{a_2-a_1+2-2g}
{\Bbb L ^{\dim\,\op{Hom}(\Cal M', V)}-1\over\Bbb L -1} d\mu_{\Cal M'}.$$
Consider the isomorphism
$$\sigma: \op{Pic}_{a_1}(X)\to \op{Pic}_{a_2+2-2g}(X), \quad \Cal M
\to \Cal M^*\otimes \Lambda^2 V \otimes\omega_X^*,$$
where $\omega_X$ is the sheaf of 1-forms. Since for  a rank 2 bundle
we have $V^*\simeq V\otimes(\Lambda^2 V)^*$, we find from the
Riemann-Roch theorem that
$$\dim\, \op{Hom}(\Cal M, V) - \dim\,\op{Hom}(\Cal M', V) = 
a_2-a_1+2-2g,$$
whenever $\Cal M' = \sigma(\Cal M)$.
 This means that the difference of the
coefficients at $z_1^{a_1}z_2^{a_2}$ in the two sides of (4.2.2) is
$$\int_{\Cal M\in\op{Pic}_{a_1}(X)} {\Bbb L^{a_2-a_1+2-2g}-1\over\Bbb L-1}
d\mu_{\Cal M} = \mu(\op{Pic}_{a_1}(X)) \cdot {\Bbb L^{a_2-a_1+2-2g}-1\over
\Bbb L-1}.$$
Our assumption that $\op{Pic}_1(X)(k)\neq\emptyset$ implies that
$\mu(\op{Pic}_{a_1}(X)) = \mu(\op{Pic}_0(X))$ and we find that the difference
between the two series in (4.2.2) is
$${\mu(\op{Pic}_0(X))\over\Bbb L -1} \sum_{a_1+a_2=\op{deg}(V)} (\Bbb L^{
a_2-a_1+2-2g}-1)z_1^{a_1}z_2^{a_2} =$$
$$= {\mu(\op{Pic}_0(X))\over\Bbb L-1} \biggl(
z_2^{\op{deg}(V)} \Bbb L^{\op{deg}(V)+2-2g} \delta
 \biggl({z_1\over \Bbb L^2 z_2}\biggr)
 - z_2^{\op{deg}(V)}\delta\biggl({z_1\over z_2}\biggl)\biggl).$$
This proves our claim.

\proclaim{(4.2.3) Corollary} Theorems 4.1.2 and 4.1.4 hold whenever $G$ has
semisimple rank 1. 
\endproclaim

\vskip .2cm

\noindent {\bf (4.3) The general case.}
We now deduce Theorems 4.1.2 and 4.1.4 for general $G$ from the case
of semisimple rank 1, following an approach which the author
learned from V. Drinfeld (this approach is also implicit in [BG]). 

Note that  unlike the functional equation (4.1.2) for $E_{G, Q}(z)$
which is to be understood in the sense of analytic continuation
(summing a series to a rational function and expanding in 
a different region), the equation (4.1.4) for $N_{G,Q}(z)$
is supposed to hold at the level of monomials (provided 
we know that $N_{G,Q}$ is indeed a Laurent polynomial). Conversely,
suppose we know that the {\it series} $N_{G,Q}(z)$ defined by
(4.1.3) satisfies the equations (4.1.4) at the level of
monomials. Then we can {\it deduce} that $N_{G,Q}(z)$
is actually a Laurent polynomial. Indeed, the support of
$N_{G,Q}$ lies in some translation of $L_+\i L$, the cone of
dominant coweights of $[G, G]$. The equations (4.1.4), if they hold at the
level of monomials, would then imply that the support of
$N_{G, Q}(z)$ lies in the intersection, over all $w\in W$, of some translations
of $w(L_+)$. But any such intersection is finite. 

Further, the group $W$ being generated by simple reflections, we
are reduced to the following fact.

\proclaim{(4.3.1) Lemma} Let $\alpha\in\Delta_{\op{sim}}$. Then the
series $N_{G, Q}(z)$ defined by (4.1.3), satisfies the equation
(4.1.4) with $w=s_\alpha$, at the level of monomials.
\endproclaim

\pf Let $P_\alpha$ be the parabolic subgroup in $G$ defined by
one negative root $(-\alpha)$ and $M_\alpha\i P_\alpha$ be the
Levi subgroup. So $M_\alpha$ has semisimple rank 1, its
maximal torus is $T$, its Borel subgroup is
$B_\alpha = B\cap M_\alpha$ and the Weyl group is $\{1, s_\alpha\}$. 
Let $P_\alpha^{\op{ab}}=P_\alpha/[P_\alpha, P_\alpha]$ and $L_\alpha=\op{Hom}(\Bbb G_m, P_\alpha^{\op{ab}})$, so we have a homomorphism
$\lambda: L\to L_\alpha$. Since $M_\alpha$ is isomorphic to
the quotient of $P_\alpha$ by its unipotent radical, we have a 
projection $\varphi: P_\alpha\to M_\alpha$. A $B_\alpha$-bundle on $X$
has degree lying in $L$ and a $P_\alpha$-bundle has degree lying in
$L_\alpha$. 

Now, a $B$-structure of degree $a$ on a $G$-bundle $Q$ gives,
in particular, a $P_\alpha$-structure of degree $\lambda(\alpha)$.
Conversely, suppose that $\pi$ is a $P_\alpha$-structure on $Q$ and let
$(Q, \pi)$ be the corresponding $P_\alpha$-bundle. Then
a $B$-structure on $Q$ refining $\pi$ is the same as a $B_\alpha$-structure
on the $M_\alpha$-bundle $\varphi_*(Q, \pi)$. Let $\Gamma_b^{P_\alpha}(Q)$,
$b\in L_\alpha$, be the scheme of $P_\alpha$-structures on $Q$ of
degree $b$. We conclude that
$$E_{G, Q}(z) = \sum_{b\in L_\alpha} \int_{\pi\in\Gamma_b^{P_\alpha}(Q)}
E_{M_\alpha, \varphi_*(Q,\pi)}(z) d\mu_\pi,$$
the equality being understood at the level of each coefficient. Therefore
$$N_{G, Q}(z)=\biggl(\prod_{\beta\in\Delta_+-\{\alpha\}} (1-\Bbb L^2 z^
{\beta^\vee})\Phi(z^{\beta^\vee})\biggr)
\sum_{b\in L_\alpha} \int_{\pi\in\Gamma_b^{P_\alpha}(Q)}
N_{M_\alpha, \varphi_*(Q, \pi)}(z) d\mu_\pi.$$
Since $N_{M_\alpha, \varphi_*(Q, \pi)}(z)$ is a Laurent polynomial,
it satisfies the functional equation (4.1.4)
with $w=s_\alpha$ at the level of monomials. As $\Delta_+-\{\alpha\}$ is permuted by 
$s_\alpha$, we finally conclude that $N_{G, Q}(z)$ satisfies (4.1.4)
with $w=s_\alpha$ at the level of monomials. Lemma is proved. 

\vfill\eject

\centerline {\bf \S 5. The Kac-Moody sheaf  and $c_2$.}

\vskip 1cm

Let  $G((\lambda))$ be the loop group of $G$. This is  an
ind-scheme over $k$ such that for any commutative $k$-algebra
 $R$ we have
$$G((\lambda))(R) = G\bigl(R((\lambda))\bigr).$$
If $R$ is a field, Garland [Ga] has constructed a central
extension of $G((\lambda))(R)$ by $R^*$ which is an algebraic
version of the minimal central extension of the loop group
of a compact Lie group [PS]. We need
the more general case when $R$ is replaced by the structure
sheaf of a smooth algebraic variety over $k$. In this
section we describe the corresponding central extension
in such a way as to make clear its relation to the second
Chern class for $G$-bundles. 

Garland's extension is induced from Matsumoto's
extension of $ G\bigl(R((\lambda))\bigr)$ by $K_2\bigl(
R((\lambda))\bigr)$ by the  tame
symbol map $K_2\bigl(R((\lambda))\bigr)\to R^*$. Similarly
to this, we use the sheaf-theoretic extension of $G$ by
$K_2$ constructed by Brylinski and Deligne [BD].

\vskip .3cm

\noindent{\bf (5.1) Generalities. }

\vskip .1cm

\noindent {\bf (5.1.1) Sheaves on the category of smooth varieties.}
Let $\Cal Sm$ be the category of smooth algebraic varieties over $k$.
We consider it as a Grothendieck site with respect to the Zariski topology.
Every scheme, or ind-scheme $Z$ gives rise to a functor (sheaf) on $\Cal Sm$
represented by $Z$,
which we will denote $\Z$. If $X$ is a smooth variety,
then the sheaf on the Zariski topology of $X$ formed by $Z$-valued
functions, will be denoted by $\Z_X$. The sheaf $\underline{A}^1$,
represented by the affine line, will be denoted by $\Cal O$.

Let $\Cal F$ be a sheaf of sets on $\Cal Sm$. A vector bundle
on $\Cal F$ is, by definition, a locally free sheaf of $\Cal O$-modules
on $\Cal F$. Thus, when $\Cal F=\Z$, a vector bundle $V$  on $\Cal F$ is
a system of data associating to any morphism $f: X\to Z$, where
$X$ is a smooth variety, a vector bundle $V_f$ on $X$ and to any pair
of morphisms $X'\buildrel g\over\rightarrow X\buildrel f\over\rightarrow Z$
(with $X, X'$ being smooth varieties), an isomorphism $g^*V_f\to V_{fg}$. These
isomorphisms are required to satisfy the obvious compatibility conditions for any
triple of morphisms $X''\buildrel g'\over\rightarrow X'\buildrel g\over\rightarrow
X\buildrel f\over\rightarrow Z$ with $X, X', X''$ being smooth varieties. 

Thus, if $Z$ is a smooth variety, then a vector  bundle on $\Z$
is the same as a vector bundle on $Z$ in the usual sense.
The same holds if $Z$ is an ind-scheme which is an inductive
limit of smooth varieties. 

\vskip .2cm

\noindent{\bf (5.1.2) $A$-groupoids.} Let $A$ be an Abelian group
and $A$-Tors be the category of $A$-torsors. This is a symmetric
monoidal category with respect to the operation $\otimes$ of tensor
product of torsors. We call an $A$-groupoid a small category
$\Cal C$ enriched in $A$-Tors, i.e., a set $\op{Ob}\,\Cal C$
together with a collection of torsors $\op{Hom}_{\Cal C}(x,y)$,
$x,y\in\op{Ob}\,\Cal C$ and composition morphisms
$$\op{Hom}_{\Cal C}(y,z) \otimes \op{Hom}_{\Cal C}(x,y)\to
\op{Hom}_{\Cal C}(x,z)$$
satisfying the usual axioms. 
\vskip .1cm

Let $S$ be a Grothendieck site, and $\Cal A$
a sheaf of Abelian groups on $S$. By  globalizing the above definition
to the topos of sheaves on $S$, we get the concept of a sheaf of $\Cal A$-groupoids.
Such a sheaf $\Cal C$ consists of a sheaf of sets $\op{Ob}(\Cal C)$ on $S$ together
with, for any two local sections $x,y\in\op{Ob}(\Cal C)(U)$, a sheaf
$\Cal Hom_{\Cal C}(x,y)$ of $\Cal A|_U$-torsors and composition morphisms
for any three local sections satisfying the usual axioms and compatible with
the restrictions. 

\vskip .1cm

\noindent {\bf (5.1.2.1) Example.} Let $Z$ be an ind-scheme.
We will call sheaves of $\Cal O^*$-groupoids on $\Cal Sm$ with the
sheaf of objects $\Z$ groupoid line bundles on $\Z$. If $Z$
is a smooth variety (or an inductive limit of such), then
a groupoid line bundle on $\Z$ is just a line bundle on $Z\times
Z$ in the usual sense, equipped with a kind of  multiplicative structure. 

\proclaim{(5.1.2.2) Proposition}
Let $\Cal B'\i\Cal B$ be sheaves of groups on $S$. Then the category of
central extensions of $\Cal B$ by $\Cal A$ trivial on $\Cal B'$, is equivalent
to the category of $\Cal B$-equivariant sheaves of $\Cal A$-groupoids
with the sheaf of objects $\Cal B/\Cal B'$.
\endproclaim

\noindent {\sl Proof:} Obvious when $S=\{pt\}$. The general case reduces to this.

\vskip .2cm

\noindent {\bf (5.1.3) $\Cal A$-gerbes.} Let $S, \Cal A$ be as before. An $\Cal A$-gerbe
is (cf. [Bre] [Bry] [Gi]) a stack $\Cal R$ of (not necessarily small)
categories on $S$ in which each sheaf $\Cal Hom_{\Cal R}(x,y)$,
$x,y\in\op{Ob}(\Cal R(U))$, is made into a sheaf of $\Cal A|_U$-torsors 
in a way compatible with restrictions and compositions. A trivial example
is the stack $\Cal A$-Tors. Any sheaf $\Cal C$ of $\Cal A$-groupoids can be embedded
into an $\Cal A$-gerbe, namely $\Cal Fun^\circ_{\Cal A} (\Cal C, \Cal A\text{-Tors})$,
the stack of contravariant functors $\Cal C\to \Cal A\text{-Tors}$ which
preserve the torsor structure on $\Cal Hom$-sheaves
We will call such functors simply $\Cal A$-functors.
 The following is a basic fact of Giraud's theory.   

\proclaim{(5.1.3.1) Proposition} (a) Equivalence classes of $\Cal A$-gerbes form
a set, naturally identified with $H^2(S, \Cal A)$. We denote
by $\gamma(\Cal R)\in H^2(S, \Cal A)$ the class of an $\Cal A$-gerbe $\Cal R$. 
The equality $\gamma(\Cal R)=0$ is equivalent to $\Cal R(S)\neq\emptyset$. 

(b) If $\Cal R, \Cal R'$ are two $\Cal A$-gerbes, then so is
$\Cal Fun(\Cal R, \Cal R')$ and $\gamma(\Cal Fun(\Cal R, \Cal R')) = 
\gamma(\Cal R')-\gamma(\Cal R)$.

\endproclaim

We also need some additional properties of this correspondence.

\proclaim{(5.1.3.2) Proposition} Let $S$ be a topological space, $X\i S$ a closed
subspace, $S^\circ = S-X$. Then:

(a) Equivalence classes of pairs $(\Cal R, x)$, where $\Cal R$ is an $\Cal A$-gerbe
and $x\in\op{Ob}(\Cal R(S^\circ))$, form a set naturally identified with
$H^2_X(S, \Cal A)$. We denote by $\gamma(\Cal R, x)\in H^2_X(S, \Cal A)$ the class of
$(\Cal R, x)$. Its image in $H^2(S, \Cal A)$ is $\gamma(\Cal R)$.

(b) Assume that the local cohomology sheaves $\Cal H^i_X(S, \Cal A)$ vanish for
$i\neq 1$ (so $H^2_X(S, \Cal A) = H^1(X, \Cal H^1_X(S, \Cal A))$). Then the
2-category formed by the $(\Cal R, x)$ as in (a), has trivial 2-morphisms, i.e.,
reduces to a usual category and this category is equivalent to the category
of $\Cal H^1_X(S, \Cal A)$-torsors. We will denote $\underline{\gamma}(\Cal R, x)$
the torsor corresponding to $(\Cal R, x)$.
\endproclaim

\noindent {\sl Proof:} Part (a) being a generality, we restrict ourselves
to giving the construction of $\underline{\gamma}(\Cal R, x)$ in (b).
As $\Cal H^2_X(S, \Cal A)=0$, the object $x$ is, locally on 
on $X$, extendable to an object defined on a neighborhood of $X$.
Let $j: S^\circ\hookrightarrow  S$ be the embedding, so 
$$\Cal H^1_X(S, \Cal A) = (j_* j^* \Cal A)/\Cal A$$
is the sheaf of ``principal parts" of sections of $\Cal A$ on $S^\circ$.
If $U\subset X$ is a small open set and $\widetilde U\subset S$ is a small
open neighborhood of $U$, then $x|_{\widetilde U-U}$ admits a lifting onto
$\widetilde U$, i.e., there is $y\in \Cal R(\widetilde U)$ such that
$\op{Hom}(y|_{\widetilde U-U}, x|_{\widetilde U-U})\neq\emptyset$. Then, this
Hom is an $\Cal A(\widetilde U-U)$-torsor. The torsor
$$\op{Hom}\bigl(y|_{\widetilde U-U}, x|_{\widetilde U-U}\bigr)\bigl/\Cal A(\widetilde U)$$
over $\Cal A(\widetilde U-U)/\Cal A(U)$ is independent, up to a canonical
isomorhism, on the choice of $y$ and this is, by definition,
$\underline{\gamma}(\Cal R, x)(U)$. The rest is left to the reader.  

\vskip .3cm

\noindent {\bf (5.2) The affine Grassmannian.} Let $G[[\lambda]]$ be the group
of Taylor loops in $G$. This is a group scheme (of infinite type)
over $k$. The affine Grassmannian of $G$ is the ind-scheme $\widehat{\op{Gr}}=
G((\lambda))/G[[\lambda]]$. In fact, for every commutative $k$-algebra $R$ we have
$$\Gr(R) = G\bigl( R((\lambda))\bigr) / G\bigl(R[[\lambda]]\bigr).$$
Further, it is known that $\Gr$ is an inductive limit of projective
algebraic varieties (the closures of the affine Schubert cells). 
Our aim in this section is to construct a central extension of sheaves
of groups on $\Cal Sm$
$$1\to \Cal O^*\to \G'\to \underline{G((\lambda))}\to 1 \leqno (5.2.1)$$
trivial on $\underline{G[[\lambda]]}$. By (5.2.1-2), this is the same as to
construct a
$\underline{G((\lambda))}$-equivariant groupoid line bundle $ C$
on $\underline {\Gr}$,
 i.e., to construct, for any smooth variety $X$ and any morphisms
$f, f': X\to \Gr$, a line bundle $C(f, f')$ and for any three morphisms
$f, f', f'': X\to \Gr$, a canonical identification $C(f, f')\otimes C(f', f'')\to
C(f, f'')$ satisfying the associativity and compatible with inverse images. We
start by recalling a geometric description of morphisms into $\Gr$. 

For a scheme $X$ let $X[[\lambda]] = S\times \op{Spf}(k[[\lambda]])$.
We will consider it as a ringed space 
$(X, \Cal O_X[[\lambda]])$. Let $X((\lambda))$ be the ringed space
$(X, \Cal O_X((\lambda)))$. Then, the following is an easy consequence of the
definitions. 

\proclaim{(5.2.2) Proposition} (a) Morphisms
 $X\to\widehat{\op{Gr}}$ are in bijection
with isomorphism classes of pairs $(P, \tau)$, where $P$ is a principal
$G$-bundle on $X[[\lambda]]$ and $\tau$ is a trivialization of
$P$ over $X((\lambda))$.

\endproclaim

\proclaim{(5.2.3) Definition} Let $X$ be a $k$-scheme.
We call a { ribbon} on $X$ a formal scheme $Y$ whose underlying ordinary
scheme is $X$ and which is locally isomorphic to $X[[\lambda]]$.
An isomorphism of ribbons is an isomorphism of formal schemes
identical on $X$.

\endproclaim

The set of isomorphism classes of ribbons on $X$
is identified with $H^1(X, \Cal Aut(X[[\lambda]]))$, 
where $\Cal Aut(X[[\lambda]])$ is the sheaf of groups on $X$ formed
by (local) automorphisms of the formal scheme
 $X[[\lambda]]$ identical on $X$. 

For a ribbon $Y$ we denote $Y^\circ = Y-X$. This is a ringed
space with the underlying space $X$ and the structure sheaf locally
isomorphic to $\Cal O_X((\lambda))$.

\vskip .1cm

Let now $S$ be a smooth algebraic variety and $X\i S$ a smooth hypersurface.
We have then a particular ribbon $Y$ on $X$, namely the formal neighborhood of
$X$ in $S$. The following is a consequence of the descent lemma of
Beauville-Laszlo [BL].

\proclaim {(5.2.4) Proposition} Let $P^\circ$ be a principal $G$-bundle on $S^\circ = 
S-X$. Then the following sets are in natural bijection:

(i) Isomorphism classes of pairs $(P, \tau)$, where $P$ os a $G$-bundle on $S$
and $\tau: P|_{S^\circ}\to P^\circ$ is an isomorphism.

(ii) Isomorphism classes of pairs $(\widehat P, \widehat\tau)$, where
$\widehat P$ is a $G$-bundle on $Y$ and $\widehat\tau: \widehat P|_{Y^\circ}\to P^\circ|_{
Y^\circ}$ is an isomorphism.

\endproclaim 

\vskip .2cm

\noindent {\bf (5.3) The relative $c_2$-bundle.} Let $S$ be a smooth
algebraic variety and $X\i S$ be a smooth hypersurface. Recall (2.1) that
we have fixed a multiplicative $\K_2$-torsor $\Phi$ on $G$ and this gives
rise to a particular central extension (2.1.1) of
$\G_S$ by $\K_{2, S}$. For a principal $G$-bundle $P$ on $S$ let
$\op{Lift}(P)$ be the $\K_{2, S}$-gerbe of liftings of the $\G_S$-torsor
$\underline{P}$ to a $\widetilde G_S$-torsor, so $\gamma(\op{Lift}(P))=
c_2(P)$. 

Let now $P, P'$ be two $G$-bundles on $S$ and $\phi: P|_{S^\circ}
\to P'|_{S^\circ}$ be an isomorphism. Then $\phi$ gives an object
of the $\K_{2, S}$-gerbe $\Cal Fun (\op{Lift}(P), \op{Lift}(P'))$,
defined over $S^\circ$. Further, the Brown-Gersten-Quillen resolution of
$\K_{2, S}$ (see [Q]) implies that
$$\Cal H^1_X(S, \K_2) = \Cal O_X^*, \quad \Cal H^i_X(S, \K_2) = 0, i\neq 1.
\leqno (5.3.1)$$
Therefore, Proposition 5.1.3.2(b) implies the following.

\proclaim{(5.3.2) Proposition}
To every $P, P', \phi$ as above, there is naturally associated a line
bundle $C(P, P', \phi)$ on $X$.The image of the class of $C(P, P', \phi)$
in $\op{Pic}(X)$ under the direct image homomorphism $\op{Pic}(X)
\to\op{CH}^2(S)$ is equal to $c_2(P')-c_2(P)$. Given three $G$-bundles
$P, P', P''$ on $S$ and isomorphisms
$P|_{S^\circ}\buildrel \phi\over\rightarrow P'|_{S^\circ}
\buildrel \phi'\over\rightarrow P''|_{S^\circ}$, there is a natural
isomorphism
$$C(P', P'', \phi')\otimes C(P, P', \phi)\to
C(P, P'', \phi'\phi),$$
and these isomorphisms are associative for any four $G$-bundles on $S$
with compatible identifications over $S^\circ$. 
\endproclaim

\noindent {\bf (5.3.3) Example.}
To a coherent sheaf $\Cal F$ on $S$ supported on (some infinitesimal neighborhood of)
$X$, there is naturally associated the determinantal line bundle $\det_X(\Cal F)$
on $X$ which is uniquely characterized by the two properties:

\vskip .1cm

(a) Multiplicativity in short exact sequences of coherent sheaves on $S$.

\vskip .1cm

(b) If $\Cal F$ is a vector bundle on $X$ of rank $r$, then $\det_X(\Cal F) =
\Lambda^r(\Cal F)$. 

\vskip .1cm

\noindent Given a vector bundle $V$ on $S$ and a locally free subsheaf $V'\i V$ coinciding 
with $V$ over $S^\circ$, we set $\det_X(V:V')=\det_X(V/V')$. By multiplicativity
one extends the construction of $\det_X(V:V')$ to the case when $V, V'$ are
two arbitrary vector bundles on $S$ identified over $S^\circ$. If now $V, V'$
have trivial determinant, so give rise to $SL_r$-bundles $P, P'$ on $S$
and to an identification $\phi: P|_{S^\circ}\to P'|_{S^\circ}$, then
$C(P, P', \phi) = \det_X(V:V')$. This statement (which can be viewed as a kind of
Riemann-Roch theorem) follows at once from Quillen's description of the
boundary map on $K_2$ in terms of ``lattices'' [Gra]. Indeed, it is this
boundary map which gives the identification $\Cal H^1_X(S, \K_2)\simeq \Cal O_X^*$. 

\vskip .2cm

Further, let $Y$ be any ribbon on a smooth algebraic variety $X$.
The Gersten conjecture for equicharacteristic regular local rings
recently proved by Panin [Pa], implies that we have, similarly to (5.3.1):
$$\Cal H^1_X(Y, \K_2) = \Cal O_X^*, \quad \Cal H^i_X(Y, K_2)=0, i\neq 1,
\leqno (5.3.4)$$
where $Y$ is considered as a scheme $\op{Spec}(\Cal O_Y)$. Because of this,
we can generalize the above construction of the relative $c_2$-bundle
$C(P, P', \phi)$ to the case when $P, P'$ are $G$-bundles over $Y$ and
$\phi$ is their identification over $Y^\circ$. Moreover, this construction
is compatible with the earlier one in the case when $Y$ is the formal
neighborhood of $X$ in $S$ and $P, P'$ come from $G$-bundles on $S$.

 \proclaim{(5.3.5) Definition}
Let $P^\circ$ be a $G$-bundle on $Y^\circ$. The $c_2$-groupoid $\frak C_2(P^\circ)$
is the sheaf of $\Cal O_X^*$-groupoids on $X$ defined as follows. An object
of $\frak C_2(P^\circ)$ over an open $U\i X$ is a pair $(P, \tau)$ where
$P$ is a $G$-bundle on $Y_U$ and $\tau: P|_{Y_U^\circ}\to P^\circ|_{Y_U^\circ}$
is an isomorphism, while $\Cal Hom_{\frak C_2(P^\circ)} ((P, \tau), (P', \tau'))$
is the $\Cal O_U^*$-torsor corresponding to the line bundle $C(P, P', \tau'\tau^{-1})$.
\endproclaim

Take now $Y=X[[\lambda]]$, so $Y^\circ = X((\lambda))$ and let $P^\circ$ be the
trivial $G$-bundle on $Y^\circ$. In this case an object of $\frak C_2(P^\circ)$
is the same a morphism $X\to\Gr$. We obtain therefore the construction
of a line bundle $C(f, f')$ on $X$ for any two morphisms $f, f': X\to\Gr$. 
The next statement is now obvious.

\proclaim{(5.3.6) Proposition} The line bundles $C(f, f')$ give rise to a groupoid
line bundle $C$ on the sheaf of sets $\underline{\Gr}\times\underline{\Gr}$ 
on $\Cal Sm$. This groupoid line
bundle is equivariant with respect to the sheaf of groups $\underline{G((\lambda))}$.
In particular, we get a central extension $\G'$ of $\underline{G((\lambda))}$ by
$\underline{\Bbb G_m}$ in the category of sheaves on $\Cal Sm$.

\endproclaim

\noindent {\bf (5.3.7) Remark.}
It is natural to expect that the central extension $G'$ can be constructed
as an group ind-scheme, not just as a sheaf on the category of smooth varieties.
For this, we need  to construct $C$ as a groupoid line bundle on the
ind-sheme $\Gr\times\Gr$ (i.e., to construct $C(f, f')$ for any
two morphisms of any scheme $X$ into $\widehat{\op{Gr}}$. 
This does not automatically follow from the above construction because
the Schubert varieties forming a direct system representing $\Gr$, are singular.
Our construction uses  the Gersten resolution and
therefore is not directly applicable to the singular case.
It is probably possible to push the construction by using the
approach of  Kumar and Mathieu [Ku]  [Ma1] (in particular, the disingularizations
of the Schubert varieties). However, it is beside the
main purpose of the present paper, where
 it is enough to consider
the sheaves of $G'$-valued functions on smooth varieties $X$ (in fact, we really
need only the case of smooth curves).  

\proclaim{(5.3.8) Definition}
Let $Y$ be a ribbon on a smooth algebraic variety $X$ and $P^\circ$ be
 a principal
$G$-bundle on $Y^\circ$. A $c_2$-theory on $P^\circ$ is an object if the
category $\Cal Fun_{\Cal O_X^*}(\frak C_2(P^\circ), \Cal O_X^*\text{-Tors})$,
 i.e., a rule $C$ which to any
local extension $P$ of $P^\circ$ to a $G$-bundle on $Y_U$, $U\i X$, associates
a line bundle $C(P)$ on $U$ and to any two such extensions $P, P'$ 
 an isomorphism
$C(P) \otimes C(P, P')\to C(P')$ in a way satisfying the associativity and
compatible with the restrictions. 
\endproclaim

\vskip .2cm

\noindent {\bf (5.4) The full Kac-Moody group.} Consider the algebraic group $\Bbb G_m$
acting on $\op{Spf}(k[[\lambda]])$ by the ``rotation of the loop'' $\lambda\mapsto z\lambda$,
$z\in \Bbb G_m$. Because of the naturality of the central extension (2.1.1) and
of the $c_2$-bundles, we have a natural action of $\underline{\Bbb G_m}$ 
on $\G'$. We define $\widehat {\G} = \G'\rtimes \underline{\Bbb G_m}$.
We regard $\widehat{\G}$ as the functor on $\Cal Sm$ represented by the
full Kac-Moody group. In particular, the maximal torus in $\widehat {\G}$
is $T_{\op{aff}}= \Bbb G_m\times T\times \Bbb G_m$ where the first $\Bbb G_m$
is the center and the second one is the group of rotations of the loop. 
Accordingly, the affine root system as introduced in (3.2) is just the
root system of $\widehat{\G}$, i.e., the system of weights of $T_{\op{aff}}$
in the adjoint representation of $\widehat{\G}$. Similarly, the action of
$\widehat W$ on $L_{\op{aff}}$ as given in (3.3.1) comes from the action
on $T_{\op{aff}}$ of the normalized of $T_{\op{aff}}$ in $\widehat{\G}$.
Let $I = \{ g(\lambda)\in G[[\lambda]]: g(0)\in B\}$ be the Iwahori subgroup
in $G((\lambda))$ and let $\underline{\widehat{B}} = \underline{\Bbb G_m}\times
\underline {I}\rtimes \underline{\Bbb G_m}\i \widehat {\G}$. Then positive affine roots
are the weights on the Lie algebra of $\underline{\widehat B}$.

More generally, let a smooth algebraic variety $X$ be fixed. 
We have then an action of the sheaf of groups $\Cal  Aut(X[[\lambda]])$
on $\G'_X$ by group automorphisms and we define $\Cal G_X = 
\G'_X\rtimes \Cal Aut (X[[\lambda]])$. We also define $\Cal B_X$
to be the semidirect product $\underline{\Bbb G_m}\times \I_X\rtimes
\Cal Aut(X[[\lambda]])$. 

\proclaim {(5.4.1) Proposition} Let $X$ be a smooth algebraic variety. Then:

(a)  The
category of sheaves of $\Cal G_X$-torsors is equivalent to the category of
triples $(Y, P^\circ, C)$ where $Y$ is a ribbon on $X$, $P^\circ$ is a $G$-bundle
on $Y^\circ$ and $C$ is a $C_2$-theory on $P^\circ$. 

(b) The category of sheaves of $\widehat{\G}_X$-torsors is equivalent to the
category of triples $(N, P^\circ, C)$ where $N$ is a line bundle on $X$,
$P^\circ$ a $G$-bundle on the punctured formal neighborhood of $X$ in $N$ and
$C$ is a $c_2$-theory on $P^\circ$. 
\endproclaim

\noindent {\sl Proof:} This is a more or less straightforward consequence of the
definitions. Thus, in (a), a $\Cal G_X$-torsor gives an $\Cal Aut (X[[\lambda]])$-torsor
via the homomorphism $\Cal G_X\to \Cal Aut (X[[\lambda]])$,
 and an $\Cal Aut(X[[\lambda]])$-torsor is the same as a ribbon, say, $Y$. Lifting
an $\Cal Aut (X[[\lambda]])$-torsor corresponding to $Y$ to a torsor over
$\underline{G((\lambda))}_X\rtimes \Cal Aut(X[[\lambda]])$ amounts to giving
a $G$-bundle $P^\circ$ on $Y^\circ$. Further lifting of a torsor over
$\underline{G((\lambda))}_X\rtimes \Cal Aut(X[[\lambda]])$ to a torsor
over $\Cal G_X$ amounts to fixing a $c_2$-theory in $P^\circ$ because 
the central extension is defined, in the first place, in terms of the relative
$c_2$-bundles. We leave the details to the reader. The next proposition is
equally straightforward. 

\proclaim{(5.4.2) Proposition} Let $\Cal Q$ be a $\Cal G_X$-torsor corresponding to
the data $(Y, P^\circ, C)$ as in (5.4.1)(a). Then,  $\Cal B_X$-structures in
$\Cal Q$ are naturally identified with isomorphism classes of
triples $(P, \tau, \pi)$ where
$P$ is a $G$-bundle on $Y$, while $\tau: P|_{Y^\circ}\to P^\circ$ is an isomorphism
and $\pi$ is a $B$-structure in $P|_X$. 
\endproclaim

Note that we have a homomorphism of sheaves of groups
$$\Cal Aut(X[[t]])\to \Cal O_X^* \leqno (5.4.3)$$
which takes an automorphism $g$ of $X[[\lambda]]$ identical on $X$, into
the invertible function given by the action of the differential of $g$
on the normal bundle of $X$ in $X[[\lambda]]$. Therefore we get
a homorphism $\Cal B_X\to \T_{\op{aff}, X}$; in particular, a $\Cal B_X$-torsor
$\Cal T$ has a degree $\op{deg}(\Cal T)$ lying in $L_{\op{aff}}= \Bbb Z\oplus L\oplus\Bbb Z$
which is the group of 1-parameter subgroups in $T_{\op{aff}}$. The above
proposition is now easily complemented as follows. 

\proclaim{(5.4.4) Proposition} Let $\Cal T$ be a $\Cal B_X$-torsor and
$(P, \tau, \pi)$ be the data corresponding to $\Cal T$ by Proposition 5.4.2.
Then
$$\op{deg}(\Cal T) = \biggl( c_2(P), \op{deg}(P|_X, \pi), X^2_Y\biggr).$$ 

\endproclaim

\vskip 2cm

\centerline {\bf \S 6. Motivic Eisenstein series for $\widehat G$.}

\vskip 1cm

\noindent {\bf (6.1) The Eisenstein series.} We now assume
that $X$ is a smooth projective curve, as in \S 2-4. Let $Y$ be a ribbon
on $X$ and $P^\circ$ be a principal $G$-bundle on $Y^\circ$.
Fix a $c_2$-theory $C$ on $P^\circ$. If $P$ is an extension of $P^\circ$
to the whole of $Y$, we will write $c_2(P) = \op{deg}\, C(P)\in \Bbb Z$.
We denote by $\Cal Q$ the $\Cal G_X$-torsor corresponding to
$(Y, P^\circ, C)$ by Proposition 5.4.1.

\proclaim{(6.1.1) Theorem} Suppose that $X^2_Y$, the self-intersection
index of $X$ in $Y$, is negative. Then for each $n\in \Bbb Z$ there exists
a scheme $\Gamma_n(\Cal Q)$ which is a fine moduli space
for extensions $P$ of $P^\circ$ to a $G$-bundle on $Y$ with $c_2(P)=n$.
The scheme $\Gamma_n(\Cal Q)$ 
is empty for $n \ll 0$.
\endproclaim

\pf This is obtained by the same arguments as in Theorem 2.2.1,
except that we should use Chern classes with values in local cohomology
and the Grothendieck-Riemann-Roch theorem for such classes. 
We leave the details to the reader. 

\proclaim {(6.1.2) Corollary} For each $n\in \Bbb Z$ and
 $a\in L$ there exists a scheme $\Gamma_{n, a}(\Cal Q)$ of finite type,
which is a fine moduli space for triples $(P, \tau, \pi)$
where $(P, \tau)$ are as in Theorem 6.1.1 and $\pi$ is a $B$-structure in
$P|_X$ of degree $a$. 
\endproclaim

The motivic Eisenstein series corresponding to $Y, \widehat G$ and $\Cal P$
is then defined to be
$$E_{\Cal G, \Cal Q}(q, z, v) = \sum_{n\in\Bbb Z, a\in L} 
\mu(\Gamma_{n, a}(\Cal Q))q^n z^a v^d, \quad d=-X^2_Y.
\leqno (6.1.3)$$
As before, we will  write $t\in T^\vee_{\op{aff}}$ for $(q,z, v)$.

\proclaim{(6.1.4) Proposition} Let $S$ be a smooth projective surface,
$X\i S$ be a smooth curve, $P^\circ$ be a $G$-bundle on $S-X$
and $Y$ be the formal neighborhood of $X$ in $S$. Fix some $c_2$-theory
$C$ in $P^\circ|_{Y^\circ}$ and let
 $\Cal Q$ be the corresponding $\Cal G_X$-torsor. 
 Then the generating function $E_{G, P^\circ}(q,z, v)$
from (3.1.1) differs from $E_{\Cal G, \Cal Q}(q,z, v)$ by a factor
of $q^m$ for some $m$.
\endproclaim

\pf  By Proposition 5.2.4, 
there is a bijection between extensions of $P^\circ$ to a bundle on $S$
and extensions of $P^\circ|_{Y^\circ}$ to a bundle on $Y$. 
So unless both $E_{G, P^\circ}$ and $E_{\widehat G, \Cal P}$ are both
identically zero, there exists an extension $P$ of $P^\circ$
to a bundle on $S$. Then, $P$  gives rise to a canonical $c_2$-theory
normalized so that its value on $P|_Y$ is $\Cal O_X$
and thus to a $\Cal G_X$-torsor $\Cal Q'$.

It follows  that $E_{\Cal G, \Cal Q'}(q,z, v) = 
q^{-c_2(P)} E_{G, P^\circ}(q, z, v)$, where $c_2(P)$ is the usual second Chern class
of $P$ on $S$. Further, $E_{\Cal G, \Cal Q}(q, z, v)$ and 
 $E_{\Cal G, \Cal Q'}(q, z, v)$
differ by a factor of a power of $q$, since the concepts of $c_2\in\Bbb Z$ for
$G$-bundles on $Y$ defined by $\Cal Q$ and $\Cal Q'$, differ by a constant,
namely, by $\op{deg}\, C(P)$. Proposition is proved. 

\vskip .1cm

In view of the above proposition, Theorems 3.4.1 and 3.4.4 would follow from
the next fact.

\proclaim{(6.1.5) Theorem} For any $Y$ with $X^2_Y=-d<0$, any $G$
and any $\Cal G_X$-torsor $\Cal Q$ as above,
the motivic Eisenstein series $E_{\Cal G, \Cal Q}(q,z, v)$
satisfies the properties claimed in 
Theorems 3.4.1 and 3.4.4.
\endproclaim

\vskip .2cm

\noindent{\bf (6.2) Proof of Theorem 6.1.5.}
Let $D(t)$ be as in (3.4.2) and $N_{\Cal G, \Cal Q}(t) = E_{\Cal G, \Cal Q}(t)
D(t)$, $t=(q,z,v)$. Note that 
$$E_{\Cal G, \Cal Q}(q, z, v) = v^d \sum_n q^n\int_{P\in \Gamma_n(\Cal Q)} E_{G, P|_X}
(z) d\mu_P,\leqno (6.2.1)$$
we find that for any $n$ the coefficient at $q^nv^d$ in $E_{\Cal G, \Cal Q}$
is a series in $z$ whose support (a subset in $L$) lies is some translation of the
cone of dominant coweights. This implies that $N_{\Cal G, \Cal Q}(t)$
is  well-defined as formal series.
Further, the desired functional equation  for $N_{\Cal G, \Cal Q}(t)$
in 
 the case $w\in W\i\widehat W$
follows from (6.2.1) and Theorem 4.1.4. So it is enough to consider the case
$w=s_{\alpha_0}$, where $\alpha_0$ is the new simple affine root. 
This is parallel to the proof of Lemma 4.3.1 except that we have to deal with
parahoric subgroups in $G((\lambda))$ instead of parabolic subgroups in $G$.

Let $\Pi_0\i G((\lambda))$ be the parahoric subgroup corresponding
to the negative affine root $(-\alpha_0)$ and $\underline{\Pi}'_X\i\G'_X$
be the preimage of $\underline{\Pi}_{0, X}\i \underline{G((\lambda))}_X$. 
Set further $\Cal P = \underline{\Pi}'_X\cdot \Cal Aut(X[[\lambda]]) \i\Cal G_X$.

Let $M_0\i\Pi_0$ be the standard Levi subgroup. We have then the projection
$\phi_0: \Pi_0\to M_0$. Let $\underline{M}'_X\i\underline{G}'_X$ be the preimage
of $\underline{M}_{0,X}$ and $\underline{M_X}=\underline {M}'_X\cdot \Cal O_X^*$,
where $\Cal O_X^*$ is regarded as the subgroup in $\Cal Aut (X[[\lambda]])$
formed by the automorphisms $\lambda\mapsto f(x)\cdot\lambda$, $f(x)\in\Cal O_X^*$. 
The projection $\phi_0$ together with the homomorphism (5.4.3) induce a
homomorphism $\phi: \Cal P\to \underline{M}_X$. 

Note that the sheaf of groups $\underline{M}_X$ is in fact formed by
regular functions with values in a reductive algebraic group $M$ over $k$,
of semisimple rank 1, with maximal torus $T_{\op{aff}}$, 
 root system $\{\pm \alpha_0\}$ and Weyl group
$\{1, s_{\alpha_0}\}$. Let $B'\i M$ be the standard Borel subgroup, so
$\underline{B}'_X= \underline{M}_X\cap \Cal B_X$. 

\proclaim{(6.2.2) Lemma}
let $\Cal Q$ be a $\Cal G_X$-torsor.
Then, defining a $\Cal B_X$-structure 
in $\Cal Q$ is equivalent to, first, defining a $\Cal P$-structure,
say, $\varpi$ and then, defining a $\underline{B}'_X$-structure
in the $\underline{M}_X$-torsor $\phi_*(\Cal Q, \varpi)$.

\endproclaim

\pf Follows from the equalities
$$\Cal P= \underline{M}_X\cdot \Cal B_X, \quad \underline {B}'_X = 
\underline{M}_X\cap\Cal B_X.$$

\vskip .1cm

Let $L'_{\op{aff}}=\op{Hom}(\Bbb G_m, M^{\op{ab}})$, so that we have a homomorphism
$L_{\op{aff}}\to L'_{\op{aff}}$ with kernel isomorphic to $\Bbb Z$. 
To finish the proof of the functional equation of $N_{\Cal G, \Cal Q}(t)$
for $w=s_{\alpha_0}$ 
 in the manner exactly identical
to the proof of Lemma 4.3.1, it is enough to establish the next
proposition. 

\proclaim{(6.2.3) Proposition} Let $b'\in \overline{L}'$. Then there
exists a scheme $\Gamma_{b'}^{\Cal P}(\Cal {Q})$ of finite type over $k$,
parametrizing $\Cal P$-structures in $\Cal Q$ of degree $b'$. 
\endproclaim

\pf Let $w\in\widehat W$ be an element taking $\alpha_0$ to a simple non-affine
root $\alpha\in\Delta_{\op{sim}}$ and $f\in G((\lambda))$ be a representative
of $w$. Then the conjugation with $f$ takes $\Pi_0$ into the parahoric subgroup
$\Pi_\alpha$ corresponding to $(-\alpha)$. More precisely,
$$\Pi_\alpha = \bigl\{g(\lambda)\in  G[[\lambda]]\,\, 
\bigl| g(0)\in P_\alpha)\bigr\},$$
where $P_\alpha\i G$ is the parabolic subgroup corresponding to 
$(-\alpha)$, see the proof of Lemma 4.3.1 from which we borrow
other notations as well. Let $\Cal P_\alpha$ be the sheaf of groups on $X$
constructed from $\Pi_\alpha$ int he same was as $\Cal P$ was constructed
from $\Pi_0$: we first lift $\underline{\Pi}_{\alpha, X}$ to
$\G'_X$ and then multiply with $\Cal Aut(X[[\lambda]])$. 
Degrees of $\Cal P_\alpha$-torsors lie
therefore in ${L}^\alpha_{\op{aff}}:= \Bbb Z\oplus L_\alpha\oplus \Bbb Z$. Since 
$\Cal P$- and $\Cal P_\alpha$-structures in $\Cal  Q$ are in bijection (via
conjugation with $f$), it is enough to show that for any
$b''=(n, b, -d)\in {L}^\alpha_{\op{aff}}$, there exists a scheme of finite
type $\Gamma_{b''}^{\Cal P_\alpha}(\Cal {Q})$ 
parametrizing $\Cal P_\alpha$-structures in $\Cal  Q$
of degree $b''$. But this is trivial: $\Gamma_{b''}^{\Cal P_\alpha}(\Cal {Q})$
is fibered over the scheme $M_{G, P^\circ}(n)$ with the fiber over $(P, \tau)$
being the scheme $\Gamma_{G, P|_X}^{P_\alpha}$ of $P_\alpha$-structures
in $P|_X$ of degree $b$.

This finishes the proof of Theorem 6.1.5 and thus of Theorems 4.1.2 and 4.1.4.

\vskip 2cm

\centerline{\bf \S 7. Explicit calculations for $X=\Bbb P^1$.}

\vskip 1cm

In this section we assume that $X=\Bbb P^1$. Our aim is to
classify all $\Cal G_X$-torsors  whose associated ribbons
have negative self-intersection index and, for each such torsor,
to find the corresponding Eisenstein series completely. 
Our method can be seen as a version  of Langlands' calculation
of the constant term of the Eisenstein series, but applied
to  the case of Kac-Moody groups and to the framework of motivic measures.

\vskip .3cm

\noindent {\bf (7.1) Grothendieck's theorem for $\Cal G_X$-torsors.}
Recall that we have the following homomorphisms of sheaves of
groups on $X=\Bbb P^1$:
$$\T_{\op{aff}, X}\hookrightarrow \Cal G_X\to \Cal Aut(X[[\lambda]])
\buildrel \text{(5.4.3)}\over\longrightarrow \Cal O_X^*.\leqno (7.1.1)$$
If $\Cal Q$ is  a $\Cal G_X$-torsor and $(Y, P^\circ, C)$ are the data
corresponding to $\Cal Q$ by Proposition 5.4.2, then
the image of the class of $\Cal Q$ under the map
$H^1(X, \Cal G_X)\to H^1(X, \Cal O_X^*)$ is the class of the normal
bundle $N_{X/Y}$. We will say that $\Cal Q$ is of negative index,
if $N_{X/Y}$ has negative degree, so $N_{X/Y}=\Cal O(-d)$, $d>0$.
Since $X=\Bbb P^1$, isomorphism classes of $T_{\op{aff}}$-bundles on
$X$ are in bijection with $L_{\op{aff}}=\op{Hom}(\Bbb G_m, T_{\op{aff}})$.
For $a\in L_{\op{aff}}$ be denote by $\Cal O^*(a)$ the corresponding
$T_{\op{aff}}$-bundle on $X$. By $\Cal O^*(a)_{\Cal G}$ we will
mean the $\Cal G_X$-torsor induced  from $\Cal O^*(a)$
by the first homomorphism in (7.1.1).

The following fact can be seen as an analog of the theorem of
Grothendieck [Gro] for Kac-Moody groups.

\proclaim{(7.1.2) Theorem} Suppose $k$ is algebraically closed.
A $\Cal G_X$-torsor $\Cal Q$ whose index is negative, 
 has the form $\Cal O^*(b)_{\Cal G}$,
where $b\in L_{\op{aff}}$ is an antidominant coweight, which is uniquely
defined by $\Cal Q$.
\endproclaim

\pf Let $(Y, P^\circ, C)$ be the  the data
corresponding to $\Cal Q$ by Proposition (5.4.2). We start
by identifying $Y$.

\proclaim {(7.1.3) Lemma}
Any ribbon $Y$ on $X=\Bbb P^1$ such that $X^2_Y$ is negative
has a linear structure, i.e., is isomorphic to the formal
neighborhood of $X$ in the total space of some line bundle $N$
on $X$ (which is identified with $N_{X/Y}$).
\endproclaim

\pf An automorphism of $X[[\lambda]]$ is the same as an algebra
automorphism of $\Cal O_X[[\lambda]]$ which is continuous in the
$\lambda$-adic topology and is identical on $\Cal O_X$. Such an
automorphism is uniquely determined by its values on $\Cal O_X$
and $\lambda$ which have the form
$$a\mapsto a+\sum_{i=1}^\infty D_i(a) \lambda^i ,
\quad \lambda\to \lambda+ \sum_{i=1}^\infty b_i\lambda_i,  
\leqno (7.1.4) $$
where $D_i: \Cal O_X\to\Cal O_X$ are morphisms of sheaves of vector
spaces and $b_i$ are functions such that  $b_1$ is invertible.

The only constraint on these data is that the first correspondence 
in (7.1.4) defines an algebra homomorphism $\Cal O_X\to\Cal O_X[[\lambda]]$.
This gives a set of Leibniz-type conditions on the $D_i$ which are sometimes
expressed by saying that $(D_i)$ forms a higher derivation of
$\Cal O_X$, see [Ma3], \S 27. In particular, $D_1$ is a derivation
(vector field) in the usual sense. Further, if $D_1=...=D_{i-1}=0$,
then $D_i$ is a section of a line bundle on $X$, namely the $i$th
tensor power of the tangent bundle $T_{\Bbb P^1}=\Cal O(2)$.
For $i\geq 1, j\geq 2$ let $F^{ij}$ be the sheaf of subgroup
in $\op{Ker}(\phi)$ defined by $D_1, ..., D_{i-1}=0$,
$b_2, ..., b_{j-1}=0$. This is a decreasing filtration 
by normal subgroups with the quotients  being Abelian and,
more precisely, $\op{gr}^{ij}_F = \Cal O(2i)$.
Now, $\Cal Aut(X[[\lambda]])$ is a semidirect
product of $\Cal O_X^*$ and $\op{Ker}(\phi)$.
The action of $\Cal O_X^*$ on  the line bundle
$\op{gr}^{ij}_F$ induced by the conjugation is
via the homomorphism $\Cal O_X^*\to\Cal O_X^*$
given by raising to the $(i+j)$th power.
Our lemma says that the preimage in $H^1(X, \Cal Aut(X[[\lambda]]))$
 of  the class of $\Cal O(-d)$ in
$H^1(X, \Cal O^*)$ consists of one element, if $d>0$. 
To see this, we use the filtration of $\Cal Aut(X[[\lambda]])$
formed by $G^{ij}= F^{ij}\rtimes\Cal O_X^*$ and find that possible liftings
from
 $\Cal Aut(X[[\lambda]]/G^{ij}$ to $\Cal Aut(X[[\lambda]])/G^{i+1,j}G^{i, j+1}$
form a homogeneous space over $H^1(X, 2i+d(i+j))=0$. 
QED.

\vskip .2cm

We denote by $Y_d$ the formal neighborhood of $X$ in the total space
of the line bundle $\Cal O(-d)$. Thus we have a projection
$p_d: Y^\circ_d\to X$. Note that we have an isomorphism
of line bundles $p_d^*\Cal O(d)\to\Cal O_{Y_d^\circ}$ on $Y_d^\circ$.

\proclaim{(7.1.5) Lemma} The Picard group of $Y_d^\circ$ is identified
with $\Bbb Z/d$ and consists of line bundles $p_d^*\Cal O(i)$, where $i$ is
taken modulo $d$.
\endproclaim

\pf  Consider first the case $d=1$. Then we have the blow-down
morphism $\sigma: Y_1\to D$ where $D=\op{Spec} \, k[[x,y]]$ is the 2-dimensional
formal disk. Let $D^\circ = D-\{0\}$ be the punctured formal disk
and $j: D^\circ\hookrightarrow D$ the embedding. Since $\sigma$ is an isomorphism
outside $X$, we find that a line bundle on $Y_1^\circ$ is the same
as a line bundle on $D^\circ$. But we have the following fact
which expresses the well known property that a reflexive sheaf
on a smooth surface is a vector bundle. 

\proclaim {(7.1.6) Lemma} If $V^\circ$ is any vector bundle on $D^\circ$,
then $j_*V^*$ is a vector bundle on $D$. In particular (since any vector
bundle on $D$ is trivial), $V^\circ$ is trivial.
\endproclaim

Thus the case $d=1$ is clear. If $d$ is arbitrary, consider the
morphism of the total space of $\Cal O(-1)$ to the total
space of $\Cal O(-d)=\Cal O(-1)^{\otimes d}$ given at each fiber by
$v\mapsto v\otimes ... \otimes v$ ($d$ times). Let $\phi$
be the induced morphism $Y_1^\circ\to Y_d^\circ$. This is a Galois
covering  with the Galois group being the group scheme $\mu_d$ of $d$th
roots of unity. If $\Cal L$ is a line bundle on $Y_d^\circ$, then
$\phi^*\Cal L$ is  trivial, so $\Cal L$ can be defined, by descent,
by  defining a $\mu_d$-action in the trivial bundle on $Y_1^\circ$.
By doing so, we get exactly the bundles $p_d^*\Cal O(i)$. QED

To finish the  proof of Theorem 7.1.2 it is enough to establish
the following fact.

\proclaim {(7.1.7) Proposition} Let $d>0$. Then the following sets
are naturally identified:

(i) Isomorphism classes of $G$-bundles on $Y_d^\circ$.

(ii) Conjugacy classes of homomorphisms $\mu_d\to G$.

(iii) Dominant  affine coweights of $G$ of the form $(0, a, d)\in L_{\op{aff}}
=\Bbb Z\oplus L \oplus \Bbb Z$. 

\endproclaim

\noindent {\bf (7.1.8) Remark.}
Note that the definition of the set (iii) can be phrased in a more
suggestive form: as  the set of isomorphism classes of level
$d$ irreducible projective representations of the loop group $G^L((\lambda))$
where $G^L$ is the Langlands dual group of $G$.
 Later in this section we will interpret the Eisenstein-Kac-Moody
series corresponding to $Y_d$ and a $G$-bundle $P^\circ$ on $Y_d^\circ$
 as a deformation of the character of the representation corresponding,
by the above, to $P^\circ$. 
A natural problem
is then to construct the representation itself in some algebro-geometric
terms. Note the similarity of our situation with the framework of
Nakajima [N]: the weight  of an irreducible representations is
in both cases encoded in the topological type of a bundle on some
open part of a variety (neighborhood of the infinity on an ALE space
in Nakajima's case, the complement of $X$ in $Y_d$ in our situation).

\vskip .1cm

\noindent {\bf Proof of (7.1.7):}
 If $d=1$, all three sets consist of one element, so the statement is true.
If $d>1$, we consider the morphism $\phi_d: Y_1^\circ\to Y_d^\circ$ introduced
in the proof of Lemma 7.1.5. Let $P^\circ$ be a $G$-bundle on $Y_d^\circ$
Then $\phi_d^*P^\circ$ is trivial, as can be seen by applying Lemma
7.1.6 to any vector bundle on $Y_d^\circ$ associated to $P^\circ$ and
a representation of $G$. Thus $P^\circ$ is obtained, by descent,
by defining a $\mu_d$-action in the trivial $G$-bundle on $Y_1^\circ$.
Isomorphism classes of such descent data form exactly the set (ii). 
Further, any homomorphism $\mu_d\to G$ is factored
through $T$. (See [W] for representation theory of $\mu_d$
in finite characteristic which is the same as representation
theory of $\Bbb Z/d$ in characteristic 0.) 
 This shows that $P^\circ$ has a $T$-structure. 
The identification of conjugacy classes
(with respect to $G$ or $W$)  of such $T$-structures, or,
what is the same, of the set (ii), with (iii),  is straightforward
(compare with Kac's classification of automorphisms of finite
order of a semisimple Lie algebra [Ka]).   

\vskip .3cm

\noindent {\bf (7.2) The affine Hall polynomials. }
Let $A$ be a field and $l\in A$ be a nonzero element.
 We 
introduce the twisted $\widehat W$-action on $T^\vee_{\op{aff}, A}$ given by
$$w*t = l^{-\widehat{\rho}} w(l^{\widehat{\rho}} \cdot t),
\leqno (7.2.1)$$
so that, for example, for $\alpha\in\widehat{\Delta}_{\op{sim}}$ we have
$s_\alpha *t  = l^{-\alpha} s_\alpha(t)$. Set
$$K(t; l) = \prod_{\alpha\in\widehat \Delta_+}{1- lt^{\alpha^\vee}\over
1-t^{\alpha^\vee}}.\leqno (7.2.2)$$
We consider this as a power series in $t=(q,z,v)\in T^\vee_{\op{aff}}=
\Bbb G_m\times T^\vee\times\Bbb G_m$. As in (3.3), we can consider
$K$ as a function on $(T^\vee\times \Bbb G_m)_{A((q))}$ and as such,
it is a meromorphic function. 

For an antidominant $b\in L_{\op{aff}}$ we define the affine Hall
polynomial to be
$$P_b(t; l) = \sum_{w\in\widehat W} w*(t^b K(t)),\quad  t\in T_{\op{aff}}.
\leqno (7.2.3)$$
This definition is entirely similar to the definition of
Hall polynomials for reductive groups as given by Macdonald [Mac]. 
Unlike the finite-dimensional theory, $P_b(t, l)$ is not a polynomial
in $t$ any more; however, the following fact holds.

\proclaim{(7.2.4) Proposition} (a) $P_b(t; l)$ is an analytic function
on the rigid analytic space $(T^\vee\times \Bbb G_m)_{
A((q))}^{\op{an}}$.

(b) $P_b(t; 1)$ is the monomial symmetric function (the 
result of averaging of
$t^b$ over $\widehat W$).

(c) If $k=\Bbb C$, then  $P_b(t; 0)= \chi_{-b}(t^{-1})$ where
$\chi_b$ is the character of the irreducible representation of
the Kac-Moody Lie algebra (associated to $G^L$) with highest weight $(-b)$. 
\endproclaim

\pf (a) follows because the denominator of $K(t; l)$ is a
 $\widehat W$-antisymmetric function. Part (b) is obvious, which (c)
is just the Weyl-Kac character formula.  

\vskip .1cm

\noindent {\bf (7.2.5) Remark.} As shown by Macdonald, see [Mac],
any finite root system gives rise to a 2-parameter family
of polynomials on the maximal torus invariant with respect to
the Weyl group. These include, as special cases,  Hall
polynomials and Jack polynomials. The extension of the full 2-variable
Macdonald theory to affine root systems has not yet been developed.
The affine analogs of Jack polynomials have been studied
by Etingof and Kirillov [EK].

For future reference we will need the following fact.

\proclaim{(7.2.6) Proposition} We have the equality
$$P_b(t; l) = K(t; l) \sum_{w\in \widehat W} l^{l(w)} \bigl( l^{-\widehat\rho}
w(l^{\widehat\rho} t)\bigr)^b \prod_{\alpha\in\widehat\Delta_+\cap w^{-1}
(\widehat\Delta_-)} {1-t^{\alpha^\vee}\over 1-l^2 t^{\alpha^\vee}}.$$
\endproclaim 

\pf Let us, for simplicity, drop $l$ from the notation for $P_b$ and $K$.
What we need to show is that for any $w\in\widehat W$ the $w$th summand in
the right hand side of the proposed equality is equal to
$K(t)^{-1} w*(t^bK(t))$. We will work out here the case $w=s_\alpha$, $\alpha\in
\widehat\Delta_{\op{sim}}$, the general case being similar.
In this case we claim that
$$l{1-t^{\alpha^\vee}\over 1-l^2 t^{\alpha^\vee}} = {s_\alpha* \kappa(t)\over
\kappa(t)}, \quad \kappa(t) = {1-lt^{\alpha^\vee}\over 1-t^{\alpha^\vee}}.$$
But this is elementary: by using the equality $\langle \alpha, \alpha^\vee
\rangle =2$,
we see that 
$$s_\alpha*\kappa(t) = {1-l(l^{-\alpha}s_\alpha t)^{\alpha^\vee}\over
1-(l^{-\alpha} s_\alpha t)^{\alpha^\vee}} = {1-l^{-1}t^{-\alpha^\vee}\over
1-l^{-2}t^{-\alpha^\vee}}=$$
$$={t^{\alpha^\vee} -l^{-1}\over t^{\alpha^\vee}-l^{-2}} = 
{1-lt^{\alpha^\vee}\over l^{-1} (1-l^2t^{\alpha^\vee})}.$$

\vskip .2cm

\noindent {\bf (7.3) The identification of the Eisenstein series.}
Let $\Cal Q$ be a $\Cal G_X$-torsor of negative index
of the form $\Cal O^*(b)_{\Cal G}$, where $b\in L_{\op{aff}}$
is  an antidominant affine coweight. When $k$ is
algebraically closed, any torsor is of this form, by Theorem 7.1.2.
 We assume fixed a field $A$ and an $A$-valued
 motivic measure $\mu$ on $\op{Sch}_k$; in particular, $\Bbb L= \mu(A^1)$.  
 The second main result of this paper is the following.

\proclaim {(7.3.1) Theorem} 
The Eisenstein series $E_{\Cal G, \Cal Q}(t)$ associated to $\Cal Q$,
is equal to $K(t; \Bbb L)^{-1} P_b(t; \Bbb L)$.
\endproclaim

\pf  By our assumptions, $b=(m, a, -d)\in L_{\op{aff}}=
\Bbb Z\oplus L\oplus\Bbb Z$, with  $d>0$. 
In other words, in the data $(Y, P^\circ, C)$ corresponding to $\Cal Q$
by Proposition 5.4.1(a), the ribbon $Y$ is $Y_d$, the formal neighbohrood
of $X=\Bbb P^1$ in the total space of $\Cal O(-d)$. Thus, by (5.4.1)(b),
the torsor  $\Cal Q$
comes from a $\widehat G_X$-torsor which we denote $Q$. 
Let $\widehat F= G((\lambda))/I$ be the affine flag variety of $G$.
It follows that we have a  fibration $\widehat{\Cal F}$ over $X$ with
fibers isomorphic to $\widehat F$ and $\Cal B$-structures in $\Cal Q$
are just sections of this fibration.

Recall (affine Bruhat decomposition) that for any field
$\frak k$,  any two $\frak k$-points $(f, f')$ on 
$\widehat{ F}$ have a uniquely defined relative position $\delta(f, f')\in 
\widehat W$. Given $f$, we define
the affine Schubert cell  $U_w(f)$ to consist of $f'$ with
$\delta(f, f')=w$. This is an algebraic variety over $\frak k$,
isomorphic to the affine space of dimension $l(w)$.

Take $\frak k = k(X)$, the field of rational functions on $X$. 
The fibration $\widehat{\Cal F}$ is trivial over $\frak k$,
 because the bundle $P^\circ$ on $Y_d^\circ$  is trivial
over the punctured formal neighborhood of $A^1$ in $Y_d$. As $\Cal Q$
comes from a $T_{\op{aff}}$-bundle on $X$, it comes equipped with
a canonical $\Cal B$-structure $\pi_0$ which we regard as
a distinguished section of $\widehat{\Cal F}$. 
Recall that the Eisenstein series associated to $\Cal Q$ is
$$E_{\Cal G, \Cal Q}(t) = \sum_{a\in L_{\op{aff}}} \mu(\Gamma_a(\Cal Q))t^a,$$
where $\Gamma_a(\Cal Q)$ is the scheme of $\Cal B$-structures in $\Cal Q$
of degree $a$. Consider the stratification
$$\Gamma_a(\Cal Q) = \coprod_{w\in\widehat W} \Gamma_a^w(\Cal Q),$$
where $\Gamma_a(\Cal Q)$ consists of those $\Cal B$-structures in $\Cal Q$
which, being regarded as sections of $\widehat{\Cal F}$ are, over the
generic point of $X$ (i.e., over $\frak k$),
 in relative position $w$ with $\pi_0$. Then we have
$$E_{\Cal G, \Cal Q}(t) = \sum_{w\in \widehat W} E^w(t), \quad 
E^w(t) := \sum_{a} \mu(\Gamma_a^w(\Cal Q))t^a.$$
We proceed to find each $E^w(t)$ separately.  

Consider a point $x\in X$ and let $\widehat{\Cal F}_x$ be the fiber
of $\widehat{\Cal F}$ over $x$. It contains a distinguished point
$\pi_{0, x}$, the value of $\pi_0$ at $x$. Therefore
we have the Schubert cell $U_w(\pi_{0, x})\i \widehat{\Cal F}_x$.

\proclaim {(7.3.2) Lemma} The Bruhat decomposition gives a canonical
identification of algebraic varieties
$$U_w(\pi_{0, x}) \longrightarrow \bigoplus_{\alpha\in\widehat\Delta_+
\cap w^{-1}(\widehat\Delta_-)} \Cal O_{P^1}(-\langle \alpha, b\rangle)_x,$$
where $\Cal O_{P^1}(-\langle \alpha, b\rangle)_x$
is the fiber at $x$ of the line bundle $\Cal O_{P^1}(-\langle \alpha, b\rangle)$.
\endproclaim

\pf An exercise in Bruhat decomposition. Left to the reader. 

\proclaim{(7.3.3) Corollary}
We have an identification of sets
$$\coprod_{a\in L_{\op{aff}}} \Gamma_a^w(\Cal Q) \longrightarrow
\prod_{\alpha\in\widehat\Delta_+
\cap w^{-1}(\widehat\Delta_-)} \Gamma_{\op{rat}}(X, \Cal O(-\langle \alpha, b
\rangle)),$$
where $\Gamma_{\op{rat}}$ denotes
 the space of rational sections of a line bundle. 

\endproclaim

The next statement can be seen as a geometric analog of the Gindikin-Karpelevic
formula which is the the crucial ingredient in
 evaluation of the constant
term of the Eisenstein series [Lan].

\proclaim{(7.3.4) Proposition} Let $f_\alpha$, $\alpha\in \widehat\Delta_+
\cap w^{-1}(\widehat\Delta_-)$, be rational sections of $\Cal O_{P^1}(-
\langle  \alpha, b\rangle)$ and let $D_\alpha$ be the divisor of poles
of $f_\alpha$. Then the degree of the $\Cal B$-structure corresponding to
the system $(f_\alpha)$ by Corollary 7.3.3, is equal to
$$w(b) + \sum_{\alpha\in\widehat\Delta_+
\cap w^{-1}(\widehat\Delta_-)} \op{deg}(D_\alpha) \cdot\alpha^\vee.$$

\endproclaim

\pf This can be viewed as  a statement about two  arbitrary
sections $\pi, \pi'$ of $\widehat {\Cal F}$
in a generic relative position $w$: it expresses how their
degrees are related if we use the identification of the Schubert cell
$U_w(\pi)$ given by the Bruhat decomposition. Take
a reduced decomposition $w=s_{\alpha_1} ... s_{\alpha_l}$,
$l=l(w)$, $\alpha_i\in\widehat \Delta_{\op{sim}}$ and choose intermediate
sections $\pi=\pi_0, \pi_1, ..., \pi_l=\pi'$ of $\widehat{\Cal F}$ so that
generically $\delta(\pi_i, \pi_{i+1}) = s_{\alpha_i}$
(this is possible because the fibration is rationally trivial).
We find that it is enough to prove the required statement
only for $(\pi_i, \pi_{i+1})$ for each $i$. In other words,
Proposition 7.3.4 in full generality follows from the particular
case  of simple reflection $w=s_\alpha$ which we now assume. This
case, however, reduces to the case of two sections of a flag
fibration corresponding to the group $GL_2$.

In other words, we need to consider a rank 2 vector bundle
$V= \Cal O(b_1) \oplus  \Cal O(b_2)$ on $X=\Bbb P^1$, and two sections
rank 1 subbundles in $V$, namely $V_0 = \Cal O(d_1)$ and
$V_1$ generically transversal to $\Cal O(d_1)$. The
identification of Corollary 7.3.3 in this case means simply that
$V_1$ is the graph of a rational morphism $\phi: \Cal O(d_2)\to
\Cal O(d_1)$ and the required particular case of Proposition 7.3.4 is:

\proclaim {(7.3.5) Lemma} If $D$ is the divisor of poles of $\phi$, then
the subbundle $V_1$  corresponding to the graph of $\phi$ has degree
 $d_2-\op{deg}(D)$.
\endproclaim

\pf We consider the projection $V_1\to \Cal O(d_2)$. This is a morphism
of line bundles with the divisor of zeroes being exactly $D$, whence
the statement. 

\proclaim{(7.3.6) Definition} If $\Cal L$ is a line  bundle on $X$ and $D'\leq D$
are positive divisors on $X$, then we denote by $\Gamma^{D'}(X, \Cal L(D))$
the open subvariety in the affine space $\Gamma(X, \Cal L(D))$
formed by those sections whose divisor of poles (in the sense of meromorphic
sections of $\Cal L$) is equal to $D'$.
We introduce the following generating function
$$\psi_{\Cal L}(u) =
 \sum_{n\geq 0} \int_{D\in X^{(n)}} \mu(\Gamma^D(X, \Cal L((D))) u^d \quad\in\quad A[[u]].$$
\endproclaim

The geometric Gindikin-Karpelevic formula implies at once:

\proclaim{(7.3.7) Corollary} We have the equality
$$E^w(t) = z^{w(b)} \prod_ {\alpha\in\widehat\Delta_+
\cap w^{-1}(\widehat\Delta_-)} \psi_{\Cal O(-\langle b, \alpha\rangle)}
(t^{\alpha^\vee}).$$
\endproclaim

In order to finish the proof of Theorem 7.3.1 it is enough, by the above
corollary and Proposition 7.2.6, to establish the following fact.

\proclaim{(7.3.8) Proposition} If $\Cal L=\Cal O_{\Bbb P^1}(m)$, $m\geq 0$,
then
$$\psi_{\Cal L}(u) = \Bbb L^{m+1} {\zeta(\Bbb Lu)\over \zeta(u)}, \quad
\text{where} \quad \zeta(u) = {1\over (1-u)(1-\Bbb Lu)}$$
is the motivic zeta-function of $\Bbb P^1$.
\endproclaim

Indeed, we apply (7.3.8) to $m=-\langle b, \alpha\rangle$ (which
is nonnegative since $b$ is antidominant) and notice that
$$ {\zeta(\Bbb Lu)\over \zeta(u)} = {1-u\over 1-\Bbb L^2 u},$$
and that  Proposition
7.2.6 identifies $P_b(k; \Bbb L)$ with the sum of
products of factors of exactly this kind. 

\vskip .1cm

\noindent {\sl Proof of (7.3.8):} 
Note that $\Gamma(X, \Cal L(D)) = \bigcup_{D'\leq D} \Gamma^{D'}
(X, \Cal L(D))$. Now
 let us write
$$\zeta(u)\psi_{\Cal L}(u) = \biggl( \sum_{n'\geq 0} u^{n'} \int_{
D'\in X^{(n')}} d\mu\biggr) \biggl( \sum_{n''\geq 0} \int_{D''\in X^{(n'')}}
\mu(\Gamma^{D''}(X, \Cal L(D'')))d\mu\biggr)$$
$$= \sum_{n\geq 0}  u^n \int_{D\in X^{(n)}} \mu (\Gamma(X, \Cal L(D)))d\mu.$$
Because $X=\Bbb P^1$ and $\Cal L=\Cal O(m)$ with $m\geq 0$, we find
that $\dim(\Gamma(X, \Cal L(D))) = m+n+1$, so the above sum is equal to
$$\sum_{n\geq 0} u^n \mu(X^{(n)}) \Bbb L^{m+n+1} = \Bbb L^{m+1}\zeta(\Bbb Lu).$$
 This finishes the proof of Theorem 7.3.1. 

\vskip .3cm

\noindent {\bf (7.4) The parahoric Eisenstein series and the
universal blowup functions.}
We start with general definitions which make sense for
any smooth projective curve $X$. Let $\Cal Q$ be a $\Cal G_X$-torsor
and $(Y, P^\circ, C)$ be the data corresponding to $\Cal Q$
by Proposition 5.4.1(a). The parahoric Eisenstein series
associated to $\Cal Q$ is the generating function
$$F_{\Cal G, \Cal Q}(q) = \sum_{n\in\Bbb Z} \mu(\Gamma_a(\Cal Q)) q^n, \leqno
(7.4.1)$$ 
where $\Gamma_a(\Cal Q)$ is the scheme from Theorem 6.1.1. 
This is a ribbon analog of the generating function $F_{G, P^\circ}(q)$
from (2.4.1). In fact, we have the following  analog of
Proposition 6.1.4, which is proved in the same way.

\proclaim{(7.4.2) Proposition} Let $S$ be a smooth projective surface,
$X\i S$ be a smooth curve, $P^\circ$ be a $G$-bundle on $S-X$
and $Y$ be the formal neighborhood of $X$ in $S$. Fix some $c_2$-theory
$C$ in $P^\circ|_{Y^\circ}$ and let
 $\Cal Q$ be the corresponding $\Cal G_X$-torsor. 
 Then the generating function $F_{G, P^\circ}(q)$
 differs from $F_{\Cal G, \Cal Q}(q)$ by a factor
of $q^m$ for some $m$.
\endproclaim

The function $F_{\Cal G, \Cal Q} (q)$ is related to 
 the Eisenstein series  $E_{\Cal G, \Cal Q}((q,z,v)$
in the same way as the parahoric subgroup $G[[\lambda]]\i G((\lambda))$
is related to the Iwahori subgroup $I$. More precisely, considering
the sheaf of subgroups $\Cal K_X = \Cal O_X^* \times \underline {G[[\lambda]]}_X
\rtimes\Cal Aut (X[[\lambda]])$ in $\Cal G_X$, we find that the schemes
$\Gamma_n(\Cal Q)$ parametrize $\Cal K_X$-structures in $\Cal Q$. 
The interest to this particular version of the Eisenstein series
is explained by the following fact.

\proclaim{(7.4.3) Proposition}
Let $X=\Bbb P^1$ and $\Cal Q$ be the $\Cal G_X$-torsor corresponding,
by Theorem 7.1.2, to the antidominant coweight $(0, 0, -1)$. Then 
 Then
$F_{\Cal G, \Cal Q}(q)$ is the universal blowup function for the
$S$-duality theory with the gauge group $G$.
\endproclaim

Let us first explain the meaning of the proposition, cf. [LQ1-2].
 Consider
a smooth projective surface $S'$, a point $p\in S'$ and the blow-up
$\sigma: S\to S'$ with the exceptional divisor $X= \sigma^{-1}(p)$. 
Let $H$ be an ample divisor on $S'$ and $M_{G}^H(S', n)$ be the moduli
space of $H$-semistable $G$-bundles on $S'$ with $c_2=n$.  For $i\in \Bbb Z$
consider the divisor
$H_i =i \sigma^*H - X$ on $S$. Then it is known that
$H_i$ is ample for  $i\gg 0$ and further, for any given $n$
the spaces  $M_G^{H_i}(S, n)$ are  identified. This space
is denoted by $M_G^{H_\infty}(S, n)$. The statement of Proposition 7.4.3
is, more precisely, that
$$\sum_{n} \mu(M_G^{H_\infty}(S, n)) q^n = F_{\Cal G, \Cal Q}(q) \cdot
\sum_n \mu (M_G^H(S', n)) q^n. \leqno (7.4.4)$$

\vskip .1cm

\noindent {\sl Proof of (7.4.3-4):} The data $(Y, P^\circ, C)$ corresponding
to our choice of $\Cal Q$ by Proposition 5.4.1(a), are as follows: 
 $Y=Y_1$ is the formal neighborhood of $X = \Bbb P^1$ in the total
space of  $\Cal O(-1)$. The $G$-bundle $P^\circ$ on $Y^\circ$ is trivial,
and  $C$ is the $c_2$-theory on $P^\circ$ normalized so that
its value on the trivial $G$-bundle on $Y$ is $\Cal O_X$.

On the other hand, identifying $X$ with the exceptional fiber $\sigma^{-1}(p)
\i S$, we identify the formal neighborhood of $X$ in $S$ with $Y=Y_1$.
If $P$ is any $G$-bundle on $S$, then its restriction to $Y^\circ$
is trivial by Proposition 7.1.7. Thus $P$ is glued from a $G$-bundle
$P'$ on $S'$ and a bundle $P''$ on $Y$ via an identification
 (trivialization)
$\tau: P''|_{Y^\circ}\to \sigma^{-1}(P')|_{Y^\circ}$. Recalling
the definition of $\Gamma_n(\Cal Q)$, we find that
$M_G^{H_\infty}(S, n)$ is a union of strata $\Sigma_{i,j}$, $i+j=n$,
such that $\Sigma_{ij}$ is a (Zariski locally trivial) fibration over
$M_G^H(S, j)$ with fiber $\Gamma_i(\Cal Q)$. On the level
of generating functions this  gives the equality (7.4.4).

\vskip .2cm

We now proceed to give an explicit formula for $F_{\Cal G, \Cal Q}(q)$
when $X=\Bbb P^1$ and $\Cal Q$
 corresponds to an arbitraty antidominant weight $b\in L_{\op{aff}}$. 
For $a\in L$ we denote
$$\widehat\Delta(a) = \bigl\{ (n, \alpha\in \Bbb Z\times\Delta\bigl|
\, 1\leq n\leq \langle\alpha, a\rangle\bigr\}\leqno (7.4.5)$$ 
and set $\lambda(a) = |\widehat\Delta(a)|$.

\proclaim{(7.4.6) Theorem} If $X=\Bbb P^1$ and
$\Cal Q$ corresponds to $b=(m, f , -d)
\in \Bbb Z\oplus L\oplus \Bbb Z$, then
$$F_{\Cal G, \Cal Q}(q) = \sum_{a\in L}  
q^{\Psi(a,f) - d\Psi(a, a)/2} \prod_{(n, \alpha)\in\widehat\Delta(a)}
\Bbb L^{-\langle f, \alpha\rangle + mn+1} {1-q^n\over 1-\Bbb L^2 q^n}.$$
In particular, the universal blowup function ($m=f=0, d=1$)
is equal to 
$$\sum_{a\in L}  q^{-\Psi(a,a)/2}\Bbb L^{ \lambda(a)} 
 \prod_{(n, \alpha)\in\widehat\Delta(a)}
 {1-q^n\over 1-\Bbb L^2 q^n}.$$
\endproclaim

In the case when $\mu$ is given by the Euler characteristic  we have
$\Bbb L=1$ and
the last formula specializes to the theta-zero-value
of Example 2.4.2.

\vskip .1cm

\pf This is proved by the same method as Theorem 7.3.1 except
that we use the decomposition of the affine Grassmannian and not
the affine flag variety,  into Schubert cells. So we just indicate the main
steps.
 Along with the affine flag fibration $\widehat{\Cal F}\to X$
we have a fibration $\widehat{\Cal Gr}\to X$ with fibers isomorphic
to $\widehat{\op{Gr}}$. The section $\pi_0$ of $\widehat{\Cal F}$
defines the decomposition of the  fiber of $\widehat{\Cal Gr}$ over any
$x\in X$, 
into Schubert cells. These cells are labelled by $\widehat W/W=L$;
the cell corresponding to $a \in L$
is identified, similarly to Lemma 7.3.2, with the fiber
at $x$ of $\bigoplus_{(n, \alpha)\in\widehat\Delta(a)} 
\Cal O(-\langle f, \alpha\rangle +mn)_x$. More precisely,
we need to consider the roots entering the root
decomposition of the nilpotent radical of the parahoric subalgebra
$\frak g[[\lambda]]$, and these are $(n, \alpha)$ such that $n\geq 1$
and $\alpha\in\Delta$. Then we need to consider those of such roots
which are taken into negative affine roots by $a$ considered
as an element of $\widehat W$. These form precisely the set
$\widehat\Delta(a)$. The parahoric analog of the geometric
Gindikin-Karpelevic formula (7.3.4) says that the degree 
(i.e., the second Chern class) 
of the $\Cal K_X$-structure corresponding to a family $(f_{(n, \alpha)})$,
$(\alpha, n)\in\widehat\Delta(a)$,  of rational sections of 
$\Cal O(-\langle f, \alpha\rangle +mn)$ 
is equal to 
$$ (a\circ b)_1 + \sum_{(n, \alpha)} \op{deg}(D_{(n, \alpha)}) \cdot n,$$
where $(a\circ b)_1$ is the first component (lying in $\Bbb Z$)
of the result of action of $a\in \widehat W$ on $b= (m, f, -d)\in L_{\op{aff}}$
and $D_{(n, \alpha)}$ is the divisor of poles of $f_{(n, \alpha)}$. 
Recalling  the rule (3.3.1)  describing the $a$-action on $L_{\op{aff}}$,
we establish the theorem.

\vskip .2cm

\noindent {\bf (7.4.7) Remark.} 
In this paper we always consider
 the generating functions for the  motivic invariants of  the uncompactified
moduli spaces.
Accordingly, in the case $G=SL_2$
our blowup function differs slightly from
the function obtained in [LQ1-2] where the Uhlenbeck and Gieseker
compactifications are used.

\vfill\eject

\centerline{\bf References}

\vskip 1cm

\noindent [A]
 M. Atiyah, Instantons in two and four dimensions, {\it Comm. Math. Phys.} 
{\bf 93} (1984), 437-451.

\s

\noindent [ADK] E. Arbarello, C. De Concini, V.G. Kac, The infinite wedge representation and
the reciprocity law for algebraic curves, Proc. Symp. Pure Math.
{\bf 49}, pt.1, p. 171-190, Amer. Math. Soc. 1989. 

\s

\noindent [BD] J.-L. Brylinski, P. Deligne, Central extensions of  reductive
groups by $\Cal K_2$, preprint of IAS, Princeton NJ, available
from P. Deligne's home page at
\hfill\break  $<$www.math.ias.edu$>$.  

\s

\noindent [BG] A. Braverman, D. Gaitsgory,
Geometric Eisenstein series, preprint \hfill\break
 math.AG/9912097.

\s

\noindent [BL] A. Beauville, Y. Laszlo, Un lemme de descente,
{\it Comptes Rendus Acad. Sci. Ser. I, Math.} {\bf 320} (1995),
335-340.

\s

\noindent [Blo] S. Bloch, $K_2$ and algebraic cycles, {\it Ann. Math.} 
{\bf 99} (1974), 349-379.

\vskip .1cm

\noindent [Bre] L. Breen, On the Classification of 2-Gerbes and
2-Stacks, {\it Ast\'erisque} {\bf 225}, Soc. Math. France, 1994.

\s

\noindent [Bry] J.-L. Brylinski, Loop Spaces, Characteristic Classes
and Geometric Quantization (Progress in Math. {\bf 107}), Birkh\"auser,
Boston, 1993.

\s

\noindent  [DL] J. Denef, F. Loeser, Germs of arcs on singular algebraic varieties
and motivic integration,
{\it Invent. Math.} {\bf 135} (1999), 201-232.

\s

\noindent [EKLV] H. Esnault, B. Kahn, M. Levine, E. Viehweg,
The Arason invariant and mod 2 algebraic cycles, {\it J. AMS},
{\bf 11} (1998), 73-118. 

\s

\noindent [EK] P. I. Etingof, A. A. Kirillov, Jr. On the affine analogue
of Jack's and Macdonald's polynomials, {\it Duke Math. J.} {\bf 78} (1995),
229-256. 

\s

\noindent [EZ] M. Eichler, D. Zagier, The theory of Jacobi forms, Birkhauser, Boston,
1985. 

\s

\noindent [Ga] H. Garland, Arithmetic theory of loop groups,
{\it Publ. Math. IHES}, {\bf 52} (1980), 5-136.

\s

\noindent [Gi] J. Giraud, Cohomologie Non-Ab\'elienne, Springer-Verlag, 
1971. 

\s

\noindent [Go]  L. G\"ottsche, Theta-functions and Hodge numbers of moduli spaces of
sheaves on rational surfaces, preprint math.AG/9808007.

\s 

\noindent [Gra] D. Grayson, Higher algebraic K-theory II, Lecture
Notes in Math. {\bf 551}, p. 217-240, Springer-Verlag, 1976.

\s

\noindent [Gro] A. Grothendieck, Sur la classification des fibr\'es
holomorphes sur le sphere de Riemann, {\it Amer. J. Math.} 
{\bf 79} (1957), 121-138. 

\s

\noindent [H1] G. Harder, Minkowskische Reduktionstheorie \"uber
Funktionenk\"orper, {\it Inv. Math.} {\bf 7} (1969), 33-54. 

\s

\noindent [H2] G. Harder, Chevalley groups over function fields and
 authomorphic forms,
{\it Ann. Math.} {\bf 100} (1974), 249-306. 

\s

\noindent [Ka] V. G. Kac, Infinite-Dimensional Lie Algebras, 3rd Ed.
 Cambridge Univ.
Press, 1990. 

\s

\noindent [Ku] S. Kumar, Demazure character formula in arbitrary Kac-Moody
setting, {\it Invent. Math. } {\bf 89} (1987), 395-423.

\s

\noindent [Kl] S. Kleiman, Les theoremes de finitude pour le foncteur de
Picard, in: Lecture Notes in Math. {\bf 225} Springer-Verlag, 1972.

\s

\noindent [Lan] R.P. Langlands, Eisenstein series, Proc. Symp. Pure Math. {\bf 9}, p. 235-252, Amer. Math. Soc. 1966.

\s

\noindent [Lau] G. Laumon, Faisceaux automorphes li\'es aux series
d'Eisenstein, in: ``Automorphic Forms, Shimura varieties and L-functions"
(L. Clozel, J.S. Milne Eds.), vol. I, p. 227-281, Academic Press,
1990.

\s

\noindent [Lo] E. Loojenga, Root systems and elliptic curves,
{\it Invent. Math.} {\bf 33} (1976), 17-32.

\s

\noindent [LQ1] W.-P. Li, Z. Qin, On blowup formulae for the S-duality conjecture
of Vafa and Witten, {\it Invent. Math.} {\bf 136} (1999), 451-482.

\s

\noindent [LQ2] W.-P. Li, Z. Qin, On blowup formulae for the S-duality conjecture
of Vafa and Witten II. The universal functions, {\it Math. Research Letters,}
{\bf 5} (1998), 439-453.

\s

\noindent [Mac] I.G. Macdonald. Symmetric Functions and Orthogonal
Polynomials (University Lecture Series {\bf 12}), Amer. Math. Soc.
Providence RI 1998.
\s

\noindent [Ma 1] O. Mathieu, Formules de Caract\`eres pour les Alg\`ebres
de Kac-Moody G\'en\'erales, {\it Ast\`erisque}, {\bf 159-160},
Soc. Math. France, 1988. 

\s

\noindent [Ma2] H. Matsumoto, Sur les sous-groupes arithm\'etiques
des groupes semisimples deploy\'es, {\it Ann. ENS} {\bf 2} (1969), 
1-62. 

\s 

\noindent [Ma3] H. Matsumura, Commutative Ring Theory, Cambridge
Univ. Press, 1989.

\s

\noindent [MY] M. Maruyama, K. Yokogawa, Moduli of parabolic stable sheaves,
{\it Math. Ann.} {\bf 293} (1992), 77-99.

\s 

\noindent [N] H. Nakajima, Instantons on ALE spaces, quiver varieties
and Kac-Moody algebras, {\it Duke Math. J.} {\bf 76} (1994), 365-416.

\s

\noindent [Pa] I.A. Panin, The equicharacteristic case
of the  Gersten conjecture,  available at 
$<$www.math.uiuc.edu/K-theory$>$ as preprint $\#$389, Feb. 2000. 
\s

\noindent [PS] A. Pressley, G. Segal, Loop Groups, Oxford, Clarendon Press,
1985.

\s

\noindent [Q] D. Quillen, Higher algebraic K-theory I, Lecture
Notes in Math. {\bf 341}, p. 85-147, Springer-Verlag, 1973. 

\s

\noindent [VW] C. Vafa, E. Witten,
 A strong coupling test of S-duality, {\it Nucl. Phys. B}, {\bf 431} (1994),
3-77.
\s 

\noindent [W] W. C. Waterhouse, Introduction to Affine Group Schemes
(Graduare Texts in Math. {\bf 66}), Springer-Verlag, 1979.

\s

\noindent [Y] K. Yoshioka, The Betti numbers of the moduli space of stable sheaves
of rank 2 on $P^2$, {\it J. reine und angew. Math.} {\bf 453} (1994),
193-220.

\enddocument